

\def\LaTeX{%
  \let\Begin\begin
  \let\End\end
  \let\salta\relax
  \let\finqui\relax
  \let\futuro\relax}

\def\UK{\def\our{our}\let\sz s}
\def\USA{\def\our{or}\let\sz z}

\UK



\LaTeX

\USA


\salta

\documentclass[twoside,12pt]{article}
\setlength{\textheight}{24cm}
\setlength{\textwidth}{16cm}
\setlength{\oddsidemargin}{2mm}
\setlength{\evensidemargin}{2mm}
\setlength{\topmargin}{-15mm}
\parskip2mm


\usepackage[usenames,dvipsnames]{color}
\usepackage{amsmath}
\usepackage{amsthm}
\usepackage{amssymb}
\usepackage[mathcal]{euscript}

%
%


\definecolor{viola}{rgb}{0.3,0,0.7}
\definecolor{ciclamino}{rgb}{0.5,0,0.5}

\def\giorgio #1{{\color{red}#1}}
\def\juerg #1{{\color{red}#1}}
\def\elvis #1{{\color{blue}#1}}
\def\gianni #1{{\color{green}#1}}

\def\giorgio #1{#1}
\def\juerg #1{#1}
\def\elvis #1{#1}
\def\gianni #1{#1}




\bibliographystyle{plain}


%

\finqui

\def\Beq{\Begin{equation}}
\def\Eeq{\End{equation}}
\def\Bsist{\Begin{eqnarray}}
\def\Esist{\End{eqnarray}}

\def\Bthm{\Begin{theorem}}
\def\Ethm{\End{theorem}}

\def\Bcor{\Begin{corollary}}
\def\Ecor{\End{corollary}}
\def\Brem{\Begin{remark}\rm}
\def\Erem{\End{remark}}

\def\Bcenter{\Begin{center}}
\def\Ecenter{\End{center}}
\let\non\nonumber




\def\step #1 \par{\medskip\noindent{\bf #1.}\quad}


\def\Lip{Lip\-schitz}
\def\Holder{H\"older}

\def\aand{\quad\hbox{and}\quad}

\def\lhs{left-hand side}
\def\rhs{right-hand side}

\def\omegalimit{$\omega$-limit}


\def\bhv{behavi\our}


\def\multibold #1{\def\arg{#1}%
  \ifx\arg\pto \let\next\relax
  \else
  \def\next{\expandafter
    \def\csname #1#1#1\endcsname{{\bf #1}}%
    \multibold}%
  \fi \next}

\def\pto{.}

\def\multical #1{\def\arg{#1}%
  \ifx\arg\pto \let\next\relax
  \else
  \def\next{\expandafter
    \def\csname cal#1\endcsname{{\cal #1}}%
    \multical}%
  \fi \next}


\def\multimathop #1 {\def\arg{#1}%
  \ifx\arg\pto \let\next\relax
  \else
  \def\next{\expandafter
    \def\csname #1\endcsname{\mathop{\rm #1}\nolimits}%
    \multimathop}%
  \fi \next}

\multibold
qwertyuiopasdfghjklzxcvbnmQWERTYUIOPASDFGHJKLZXCVBNM.

\multical
QWERTYUIOPASDFGHJKLZXCVBNM.

\multimathop
diag dist div dom mean meas sign supp .


\def\accorpa #1#2{\eqref{#1}--\eqref{#2}}
\def\Accorpa #1#2 #3 {\gdef #1{\eqref{#2}--\eqref{#3}}%
  \wlog{}\wlog{\string #1 -> #2 - #3}\wlog{}}


\def\separa{\noalign{\allowbreak}}

\def\graffe #1{\mathopen\{#1\mathclose\}}

\def\<#1>{\mathopen\langle #1\mathclose\rangle}
\def\norma #1{\mathopen \| #1\mathclose \|}

\def\[#1]{\mathopen\langle\!\langle #1\mathclose\rangle\!\rangle}

\def\iot {\int_0^t}

\def\intQt{\int_{Q_t}}
\def\intQ{\int_{\QT}}
\def\iO{\int_\Omega}
\def\iG{\int_\Gamma}
\def\intS{\int_{\Sigma_T}}
\def\intSt{\int_{\Sigma_t}}
\def\intQi{\int_{\Qi}}
\def\intSi{\int_{\Si}}

\def\dt{\partial_t}
\def\dn{\partial_\nu}

\def\cpto{\,\cdot\,}

\def\checkmmode #1{\relax\ifmmode\hbox{#1}\else{#1}\fi}
\def\aeO{\checkmmode{a.e.\ in~$\Omega$}}
\def\aeQ{\checkmmode{a.e.\ in~$\QT$}}
\def\aeG{\checkmmode{a.e.\ on~$\Gamma$}}
\def\aeS{\checkmmode{a.e.\ on~$\ST$}}
\def\aet{\checkmmode{a.e.\ in~$(0,T)$}}

\def\aat{\checkmmode{for a.a.~$t\in(0,T)$}}

\def\Aet{\checkmmode{a.e.\ in~$(0,+\infty)$}}
\def\Aat{\checkmmode{for a.a.~$t\in(0,+\infty)$}}

\def\limn{\lim_{n\to\infty}}


\def\erre{{\mathbb{R}}}

\def\enne{{\mathbb{N}}}




\def\genspazio #1#2#3#4#5{#1^{#2}(#5,#4;#3)}
\def\spazio #1#2#3{\genspazio {#1}{#2}{#3}T0}

\def\L {\spazio L}
\def\H {\spazio H}

\def\C #1#2{C^{#1}([0,T];#2)}
\def\spazioinf #1#2#3{\genspazio {#1}{#2}{#3}{+\infty}0}
\def\LL {\spazioinf L}
\def\HH {\spazioinf H}


\def\Lx #1{L^{#1}(\Omega)}
\def\Hx #1{H^{#1}(\Omega)}

\def\LxG #1{L^{#1}(\Gamma)}
\def\HxG #1{H^{#1}(\Gamma)}

\def\LQ #1{L^{#1}(Q)}

\def\Luno{\Lx 1}
\def\Ldue{\Lx 2}

\def\Huno{\Hx 1}
\def\Hdue{\Hx 2}

\def\HunoG{\HxG 1}
\def\HdueG{\HxG 2}

\def\LunoG{\LxG 1}
\def\LdueG{\LxG 2}


\def\LQ #1{L^{#1}(Q)}


\let\theta\vartheta
\let\eps\varepsilon
\let\phi\varphi

\let\hat\widehat
\let\tilde\widetilde

\let\TeXchi\chi                         
\newbox\chibox
\setbox0 \hbox{\mathsurround0pt $\TeXchi$}
\setbox\chibox \hbox{\raise\dp0 \box 0 }
\def\chi{\copy\chibox}


\def\QED{\hfill $\square$}


\def\suG{_{|\Gamma}}

\def\VG{V_\Gamma}
\def\HG{H_\Gamma}
\def\WG{W_\Gamma}
\def\nablaG{\nabla_\Gamma}
\def\DeltaG{\Delta_\Gamma}
\def\muG{\mu_\Gamma}
\def\rhoG{\rho_\Gamma}
\def\tauO{\tau_\Omega}
\def\tauG{\tau_\Gamma}
\def\fG{f_\Gamma}
\def\vG{v_\Gamma}
\def\wG{w_\Gamma}
\def\gG{g_\Gamma}
\def\zetaG{\zeta_\Gamma}
\def\xiG{\xi_\Gamma}
\def\jsoluz{((\mu,\muG),(\rho,\rhoG),(\zeta,\zetaG))}
\def\soluz{(\mu,\muG,\rho,\rhoG,\zeta,\zetaG)}
\def\soluzn{(\mun,\muGn,\rhon,\rhoGn,\zetan,\zetaGn)}
\def\soluzi{(\mui,\muGi,\rhoi,\rhoGi,\zetai,\zetaGi)}
\def\soluzs{(\rhos,\rhoGs,\zetas,\zetaGs)}
\def\Mu{(\mu,\muG)}
\def\Rho{(\rho,\rhoG)}
\def\Zeta{(\zeta,\zetaG)}
\def\Xi{(\xi,\xiG)}

\def\gstar{g^*}

\def\rhoz{\rho_0}
\def\rhoGz{{\rhoz}\suG}

\def\mz{m_0}

\def\Beta{\hat\beta}
\def\BetaG{\Beta_\Gamma}
\def\betaG{\beta_\Gamma}
\def\betaeps{\beta_\eps}
\def\betaGeps{\beta_{\Gamma\!,\,\eta\eps}}

\def\betaz{\beta^\circ}
\def\betaGz{\betaG^\circ}
\def\Pi{\hat\pi}
\def\PiG{\Pi_\Gamma}
\def\piG{\pi_\Gamma}

\def\mun{\mu^n}
\def\muGn{\mu_\Gamma^n}
\def\rhon{\rho^n}
\def\rhoGn{\rho_\Gamma^n}
\def\zetan{\zeta^n}
\def\zetaGn{\zeta_\Gamma^n}
\def\Mun{(\mun,\muGn)}
\def\Rhon{(\rhon,\rhoGn)}
\def\Zetan{(\zetan,\zetaGn)}
\def\un{u^n}

\def\rhoo{\rho^\omega}
\def\rhoGo{\rho_\Gamma^\omega}
\def\Rhoo{(\rhoo,\rhoGo)}

\def\mus{\mu^s}
\def\rhos{\rho^s}
\def\rhoGs{\rho_\Gamma^s}
\def\zetas{\zeta^s}
\def\zetaGs{\zeta_\Gamma^s}
\def\Rhos{(\rhos,\rhoGs)}
\def\Zetas{(\zetas,\zetaGs)}

\def\mut{\mu^\tau}
\def\muGt{\mu_\Gamma^\tau}
\def\rhot{\rho^\tau}
\def\rhoGt{\rho_\Gamma^\tau}
\def\Mut{(\mut,\muGt)}
\def\Rhot{(\rhot,\rhoGt)}

\def\pt{p^\tau}
\def\pGt{p_\Gamma^\tau}
\def\qt{q^\tau}
\def\qGt{q_\Gamma^\tau}
\def\psit{\psi^\tau}
\def\psiGt{\psi_\Gamma^\tau}

\def\mutn{\mu^{\tau_n}}
\def\muGtn{\mu_\Gamma^{\tau_n}}
\def\rhotn{\rho^{\tau_n}}
\def\rhoGtn{\rho_\Gamma^{\tau_n}}
\def\Mutn{(\mutn,\muGtn)}
\def\Rhotn{(\rhotn,\rhoGtn)}
\def\utn{u^{\tau_n}}
\def\ptn{p^{\tau_n}}
\def\pGtn{p_\Gamma^{\tau_n}}
\def\qtn{q^{\tau_n}}
\def\qGtn{q_\Gamma^{\tau_n}}

\def\rhotnk{\rho^{\tau_{n_k}}}
\def\rhoGtnk{\rho_\Gamma^{\tau_{n_k}}}

\def\utnk{u^{\tau_{n_k}}}

\def\mui{\mu^\infty}
\def\muGi{\mu_\Gamma^\infty}
\def\rhoi{\rho^\infty}
\def\rhoGi{\rho_\Gamma^\infty}
\def\zetai{\zeta^\infty}
\def\zetaGi{\zeta_\Gamma^\infty}
\def\Mui{(\mui,\muGi)}
\def\Rhoi{(\rhoi,\rhoGi)}
\def\Zetai{(\zetai,\zetaGi)}

\def\musi{\mu_\infty}

\def\calVz{\calV_0}
\def\calHz{\calH_0}
\def\calVsz{\calV_{*0}}
\def\calVzp{\calV_0^{\,*}}
\def\calVp{\calV^{\,*}}
\def\calNO{\calN_\Omega}
\def\calNG{\calN_\Gamma}

\def\normaV #1{\norma{#1}_V}
\def\normaH #1{\norma{#1}_H}
\def\normaW #1{\norma{#1}_W}
\def\normaVG #1{\norma{#1}_{\VG}}
\def\normaHG #1{\norma{#1}_{\HG}}

\def\normaHH #1{\norma{#1}_{\calH}}
\def\normaVV #1{\norma{#1}_{\calV}}
\def\normaWW #1{\norma{#1}_{\calW}}

\def\Qi{Q_\infty}
\def\Si{\Sigma_\infty}
\def\QT{Q_T}
\def\ST{\Sigma_T}
\def\tn{t_n}

\def\hmu{\hat\mu}
\def\hrho{\hat\rho}
\def\hmuQ{\hmu_Q}
\def\hmuS{\hmu_\Sigma}
\def\hrhoQ{\hrho_Q}
\def\hrhoS{\hrho_\Sigma}
\def\hrhoO{\hrho_\Omega}
\def\hrhoG{\hrho_\Gamma}

\def\Uad{\calU_{ad}}

\def\pG{p_\Gamma}
\def\qG{q_\Gamma}

\def\mub{\overline\mu}
\def\muGb{\overline\mu_\Gamma}
\def\rhob{\overline\rho}
\def\rhoGb{\overline\rho_\Gamma}

\def\psiG{\psi_\Gamma}

\def\ub{\overline u}
\def\ubt{\overline u^{\,\tau}}

\def\longtime{longtime}
\def\Longtime{Longtime}
\def\CO{C_\Omega}
\def\CT{c_T}

\def\hm{\hat m}
\def\hR{\hat R}

\Begin{document}


%
\title{Asymptotic limits and optimal control for\\ the Cahn--Hilliard system 
  with convection\\ and dynamic boundary conditions}
\author{}
\date{}
\maketitle
\Bcenter
\vskip-1cm
{\large\sc Gianni Gilardi$^{(1)}$}\\
{\normalsize e-mail: {\tt gianni.gilardi@unipv.it}}\\[.25cm]
{\large\sc J\"urgen Sprekels$^{(2)}$}\\
{\normalsize e-mail: {\tt sprekels@wias-berlin.de}}\\[.45cm]
$^{(1)}$
{\small Dipartimento di Matematica ``F. Casorati'', Universit\`a di Pavia}\\
{\small and Research Associate at the IMATI -- C.N.R. Pavia}\\
{\small via Ferrata 5, 27100 Pavia, Italy}\\[.2cm]
$^{(2)}$
{\small Department of Mathematics}\\
{\small Humboldt-Universit\"at zu Berlin}\\
{\small Unter den Linden 6, 10099 Berlin, Germany}\\[2mm]
{\small and}\\[2mm]
{\small Weierstrass Institute for Applied Analysis and Stochastics}\\
{\small Mohrenstrasse 39, 10117 Berlin, Germany}\\
[1cm]
{\it Dedicated to our friend Prof. Dr. Pierluigi Colli\\[.1cm]
on the occasion of his 60th birthday\\[.1cm]
with best wishes}
\Ecenter
\Begin{abstract}\noindent
\juerg{In this paper, we study initial-boundary value problems for the Cahn--Hilliard
system with convection and nonconvex potential, where dynamic boundary conditions 
are assumed for both the associated order parameter and the corresponding chemical potential. While recent works addressed the case of viscous Cahn--Hilliard systems, the `pure' 
nonviscous case is investigated here. In its first part, the paper deals with the
asymptotic behavior of  the solutions as time approaches infinity. It is shown that the 
$\omega$-limit of any trajectory can be characterized in terms of stationary solutions,
provided the initial data are sufficiently smooth. The second part of the paper deals 
with the optimal control of the system by the fluid velocity. Results concerning
existence and first-order necessary optimality conditions are proved. Here, we have to
restrict ourselves to the case of everywhere defined smooth potentials. In both
parts of the paper, we start from corresponding known results for the viscous case,
derive sufficiently strong estimates that are uniform with respect
to the (positive) viscosity parameter,
and then let the viscosity tend to zero to establish the sought results for the
nonviscous case.} 
\End{abstract}
\salta
\pagestyle{myheadings}
\newcommand\testopari{\sc Gilardi \ --- \ Sprekels}
\newcommand\testodispari{\sc Asymptotic limits and control of convective   Cahn--Hilliard systems}
\markboth{\testodispari}{\testopari}
\finqui
%

\section{Introduction}
\label{Intro}
\setcounter{equation}{0}

In the recent paper \cite{CGS13},
the following initial-boundary value problem
for the Cahn--Hilliard system with convection was studied,
\Beq
  \dt\rho + \nabla\rho \cdot u - \Delta\mu = 0
  \aand
  \tauO \dt\rho - \Delta\rho + f'(\rho) = \mu
  \quad \hbox{in $\QT:=\Omega\times(0,T)$},
  \label{Isystem}
\Eeq
\juerg{where the unknowns~$\rho$ and~$\mu$ represent the order parameter and the chemical potential,
respectively, in a phase separation process \elvis{taking place in an incompressible fluid 
contained in a container $\,\Omega\subset\erre^3$}.}
In the above equations, 
$\tauO$~is a nonnegative constant,
$f'$~is the derivative of a double-well potential~$f$,
and $u$ represents \elvis{the (given) fluid velocity, which is assumed to satisfy $\,{\rm div}\,u=0\,$ in the bulk
and $\,u\cdot\nu=0\,$ on the boundary, where $\,\nu\,$ denotes the outward unit normal to
the boundary $\Gamma:=\partial\Omega$}.
Typical and physically significant examples of $f$ 
are the so-called {\em classical regular potential}, the {\em logarithmic double-well potential\/},
and the {\em double obstacle potential\/}\juerg{, which are given, in this order,} by
\Bsist
  && f_{reg}(r) := \frac 14 \, (r^2-1)^2 \,,
  \quad r \in \erre, 
  \label{regpot}
  \\
  && f_{log}(r) := \bigl( (1+r)\ln (1+r)+(1-r)\ln (1-r) \bigr) - c_1 r^2 \,,
  \quad r \in (-1,1),
  \label{logpot}
  \\[1mm]
  && f_{2obs}(r) := - c_2 r^2 
  \quad \hbox{if $|r|\leq1$}
  \aand
  f_{2obs}(r) := +\infty
  \quad \hbox{if $|r|>1$}.
  \label{obspot}
\Esist
Here, the constants $c_i$ in \eqref{logpot} and \eqref{obspot} satisfy
$c_1>1$ and $c_2>0$, so that $f_{log}$ and $f_{2obs}$ are nonconvex.
In cases like \eqref{obspot}, one has to split $f$ into a nondifferentiable convex part 
(the~indicator function of $[-1,1]$ in the present example) and a smooth perturbation.
Accordingly, one has to replace the derivative of the convex part
by the subdifferential and interpret the second identity in \eqref{Isystem} as a differential inclusion.

As far as the conditions on the boundary $\,\Gamma\,$ are concerned, 
instead of the classical homogeneous Neumann boundary conditions, 
the dynamic boundary condition for both $\mu$ and~$\rho$ were considered, namely,
\Bsist
  && \dt\rhoG + \dn\mu - \DeltaG\muG = 0
  \aand
  \tauG \dt\rhoG + \dn\rho - \DeltaG\rhoG + \fG'(\rhoG) = \muG
  \non
  \\
  &&\mbox{on \,$\ST:=\Gamma\times(0,T)$},
  \label{IdynBC}
\Esist
where $\muG$ and $\rhoG$ are the traces \juerg{on $\ST$} of $\mu$ and~$\rho$, respectively;
moreover, $\dn$~and $\DeltaG$ denote the outward normal derivative
and the Laplace--Beltrami operator on~$\Gamma$,
$\tauG$~is a nonnegative constant,
and $\fG'$ is the derivative of another potential~$\fG$.

\elvis{The associated total free energy of the phase separation process is  
the sum of a bulk and a surface contribution and has the form
\begin{align}\label{ftot}
&{\cal F}_{\rm tot}[\mu(t),\muG(t),\rho(t),\rhoG(t)] \nonumber\\[1mm] 
&:=\iO \Big(f(\rho(x,t))\,+\,\frac 12\,|\nabla\rho(x,t)|^2\,-\,\mu(x,t)\rho(x,t)
\Big)\,dx\nonumber\\[1mm] 
&+\iG\Big(\fG(\rhoG(x,t))\,+\,\frac 12\,|\nabla_\Gamma\rhoG(x,t)|^2\,-\,\muG(x,t)\rhoG(x,t)
\Big)\,d\Gamma\,,
\end{align}
for $t\in [0,T]$. Moreover, we remark that the Cahn--Hilliard type system \eqref{IdynBC} indicates
that on the boundary $\,\Gamma\,$ another phase separation process is occurring that is
coupled to the one taking place in the bulk. It is worth
 noting that the total mass of the order parameter is conserved during the separation process;
 indeed, integrating the first identity in \eqref{Isystem} for fixed $t\in (0,T]$ over $\Omega$, using the fact that
 $\,{\rm div}\,u=0\,$ in $\Omega$ and $\,u\cdot\nu=0$ on $\Gamma$, and invoking the first of the boundary conditions \eqref{IdynBC}, we readily find that
 \begin{equation}
 \label{conserve}
 \dt\Big(\iO\rho(t)+\iG\rhoG(t)\Big)=0\,.
 \end{equation}  } 

The quoted paper \cite{CGS13} was devoted to the study
of the initial-boundary value problem obtained by complementing
\eqref{Isystem} and \eqref{IdynBC} with the initial condition $\rho(0)=\rhoz$,
where $\rhoz$ is a given function on~$\Omega$.
By just assuming that the viscosity coefficients $\tauO$ and $\tauG$ are nonnegative
and that the potentials fulfill suitable assumptions and compatibility conditions,
well-posedness and regularity results were established.
Moreover, in~\cite{CGS17}, 
the study of the \longtime\ \bhv\ was addressed
in the viscous case, i.e., if $\tauO>0$ and $\tauG>0$.
More precisely, the \omegalimit\ (in~a suitable topology)
of any trajectory~$\Rho$ was characterized in terms of stationary solutions.
Finally, in~\cite{CGS14}, a~control problem was studied, the control being the velocity field~$u$,
and first-order necessary optimality conditions were derived.
Also this study was done in the viscous case.

In the present paper, we extend these results 
to the pure Cahn--Hilliard system, i.e., to the case when $\tauO=\tauG=0$.
For the \longtime\ \bhv, 
this will be done by considering only trajectories that start
from smoother initial data, and by assuming some additional summability on the velocity field.
For the control problem, \juerg{we  assume  that the potentials are of the everywhere 
regular} type,
establish the existence of optimal controls,
and derive first-order necessary optimality conditions for a given control to be optimal
in terms of a proper adjoint problem and a variational inequality.
To this end, we start from the results obtained in \cite{CGS14} by taking
$\tauO=\tauG=\tau>0$, and then let $\tau$ tend to zero.

Concerning the asymptotic behavior of Cahn--Hilliard type equations,
one can find a number of results in several directions.
From one side, the study can be either addressed 
on the convergence of the trajectories
or devoted to the existence of attractors.
On the other hand, the literature contains studies on
several variants of the Cahn--Hilliard equations, which are obtained,
e.g., by the coupling with other equations, 
like heat type equations or fluid dynamics equations,
or the changing or adding further terms, like viscosity or memory terms,
or the replacing the classical Neumann boundary conditions by other ones,
mainly the dynamic boundary conditions in the last years.
Without any claim of completeness, 
by starting from~\cite{Zheng},
we can quote, e.g.,
\cite{AbWi, ChFaPr, CGLN, CGPS2,  GiMiSchi, GiRo, GrPeSch, JiWuZh, PrVeZa, WaWu, WuZh} 
for the first type of issues,
and 
\cite{EfMiZe, MiZe1, MiZe2, EfGaZe, LiZho, Gal1, Gal2, GaGr1, GaGr2, GrPeSch, Mir1, Mir2, Mir3, Seg}
for the existence of global or exponential attractors.

As far as optimal control problems are concerned, there exist numerous recent contributions to general viscous
Cahn--Hilliard systems. In this connection, we refer to the papers 
\cite{CGPS1,CGS1,CGS11,CGS12,CGRS} for the case of standard boundary conditions and to
\cite{CFGS1,CFGS2,CGS8,CGS9,CGS10,CS,Fukao} for the case of dynamic boundary conditions.
For the `pure' case $\tauO=\tauG=0$, there are the works of \cite{WN,HW1,JZheng1,ZW}. Moreover, we refer to the papers
\cite{ZL1,ZL2,RS} that address convective Cahn--Hilliard systems, where in the latter contribution (see also its
generalization \cite{FRS} to nonlocal two-dimensional Cahn--Hilliard/Navier--Stokes systems) the fluid
velocity was chosen as the control parameter for the first time in such systems. In connection with
Cahn--Hilliard/Navier--Stokes systems, we also mention the works \cite{HHKK,HiKeWe,HW2,HW3,Medjo}.
In this paper, we study the control by the velocity of convective local Cahn--Hilliard systems with dynamic 
boundary conditions. To the authors' best knowledge, this problem 
has never been investigated before.

The paper is organized as follows. 
In the next section, we list our assumptions and notations,
recall the properties already known, 
and state our results on the \longtime\ \bhv\ and the control problem.
The corresponding proofs will be given in Sections~\ref{LONGTIME} and~\ref{CONTROL}, respectively.


\section{Statement of the problem and results}
\label{STATEMENT}
\setcounter{equation}{0}

In this section, we state precise assumptions and notations and present our results.
First of all, let the set $\Omega\subset\erre^3$ be bounded, connected and smooth.
As in the introduction, $\nu$~is the outward unit normal vector field on $\Gamma:=\partial\Omega$, 
and $\dn$ and $\DeltaG$ stand for the corresponding (outward) normal derivative
and the Laplace--Beltrami operator, respectively.
Furthermore, we denote by $\nablaG$ the surface gradient
and write $|\Omega|$ and $|\Gamma|$ 
for the volume of $\Omega$ and the area of~$\Gamma$, respectively.
Moreover, we widely use the notations
\Beq
  Q_t := \Omega \times (0,t)
  \aand
  \Sigma_t := \Gamma \times (0,t)
  \quad \hbox{for $0<t\leq+\infty$}.
  \label{defQS}
\Eeq
Next, if $X$ is a Banach space, $\norma\cpto_X$ denotes both its norm and the norm of~$X^3$, 
and the symbols $X^*$ and $\<\cpto,\cpto>_X$ stand for the dual space of $X$
and the duality pairing between $X^*$ and~$X$.
The only exception from the convention for the norms is given
by the Lebesgue spaces~$L^p$, for $1\leq p\leq\infty$, 
whose norms are denoted by~$\norma\cpto_p$. 
Furthermore, we~put
\Bsist
  && H := \Ldue \,, \quad  
  V := \Huno 
  \aand
  W := \Hdue,
  \label{defspaziO}
  \\
  && \HG := \LdueG \,, \quad 
  \VG := \HunoG 
  \aand
  \WG := \HdueG,
  \label{defspaziG}
  \\
  && \calH := H \times \HG \,, \quad
  \calV := \graffe{(v,\vG) \in V \times \VG : \ \vG = v\suG},
  \non
  \\
  && \aand
  \calW := \bigl( W \times \WG \bigr) \cap \calV \,.
  \label{defspaziprod}
\Esist
In the following, we will work in the framework of the Hilbert triplet
$(\calV,\calH,\calVp)$.
We thus have 
$$\<(g,\gG),(v,\vG)>_{\calV}=\iO gv+\iG\gG\vG
\quad\mbox{for every $(g,\gG)\in\calH$ and $(v,\vG)\in\calV$}.$$

Next, we list our assumptions. \juerg{For our first result,
we postulate for the structure of the system the following properties 
(which are slightly stronger than the analogous ones in~\cite{CGS13})}:
\Bsist
  && \hbox{$\tauO$ and $\tauG$ \ are nonnegative real numbers.}
  \label{hptau}
  \\[1mm]
  && \Beta,\, \BetaG : \erre \to [0,+\infty]
  \quad \hbox{are convex, proper and l.s.c.\ with} \quad
  \Beta(0) = \BetaG(0) = 0.
  \qquad
  \label{hpBeta}
  \\[1mm]
  \separa
  && \Pi,\, \PiG : \erre \to \erre
  \quad \hbox{are of class $C^2$ with \Lip\ continuous first derivatives.}
  \qquad
  \label{hpPi}
\Esist
Moreover, we set
\Beq
  f: = \Beta + \Pi
  \aand 
  \fG: = \BetaG + \PiG
  \label{defpot}
\Eeq
and assume that
\Beq 
  \hbox{the functions $f:=\Beta+\Pi$ and $\fG:=\BetaG+\PiG$ are bounded from below}.
  \label{hpbelow}
\Eeq
Thus, \juerg{compared with \cite{CGS13}, we have the additional assumption \eqref{hpbelow}}.
However, we remark that this assumption is fulfilled
by all of the potentials \eqref{regpot}--\eqref{obspot}.
We set, for convenience,
\Beq
  \beta := \partial\Beta \,, \quad
  \betaG := \partial\BetaG \,, \quad
  \pi := \Pi'
  \aand
  \piG := \PiG',
  \label{defbetapi}  
\Eeq
and assume that, with some positive constants $C$ and $\eta$,
\Beq
  D(\betaG) \subseteq D(\beta)
  \aand
  |\betaz(r)| \leq \eta |\betaGz(r)| + C
  \quad \hbox{for every $r\in D(\betaG)$}
  \label{hpCC}
\Eeq
as in \cite{CaCo}.
Finally, we assume that
\Beq
  \hbox{at least one of the operators $\beta$ and $\betaG$ is single-valued}.
  \label{perunicita}
\Eeq
\Accorpa\HPstruttura hptau perunicita
In \eqref{hpCC}, the symbols $D(\beta)$ and $D(\betaG)$ 
denote the domains of $\beta$ and~$\betaG$, repectively.
More generally, we use the notation $D(\calG)$ 
for every maximal monotone graph $\calG$ in $\erre\times\erre$,
as well as for the \juerg{corresponding} maximal monotone operators induced on $L^2$ spaces.
Moreover, for $r\in D(\calG)$,
$\calG^\circ(r)$ stands for the element of $\calG(r)$ having minimum modulus,
and $\calG_\eps$ denotes the Yosida regularization of $\calG$ at the level~$\eps$
(see, e.g., \cite[p.~28]{Brezis}).

For the data, we make the following assumptions:
\Bsist
  && u \in (\HH1{\Lx3})^3.
  \label{hpu}
  \\
  && \div u = 0 \quad \hbox{in $\Qi$}
  \aand
  u \cdot \nu = 0 \quad \hbox{on $\Si$}.
  \label{hpubis}
  \\
  && (\rhoz \,,\, \rhoGz) \in \calW \,, \quad
  \betaz(\rhoz) \in H
  \aand
  \betaGz(\rhoGz) \in \HG .
  \qquad
  \label{hprhoz}
  \\
  && \mz := \frac { \iO\rhoz + \iG\rhoGz } { |\Omega| + |\Gamma| }
  \quad \hbox{belongs to the interior of $D(\betaG)$} .
  \label{hpmz}
\Esist
\Accorpa\HPdati hpu hpmz
Notice that \eqref{hprhoz} implies that  $f(\rhoz)\in\Luno$ and $\fG(\rhoGz)\in\LunoG$.
Moreover, \eqref{hpu} entails that $u\in\LL\infty{\Lx3}$.

Let us come to our notion of solution,
which requires a low regularity level. 
A~solution on $(0,T)$ is a triple of pairs $\jsoluz$ that satisfies
\Bsist
  && \Mu \in \L2\calV ,
  \label{regmu}
  \\
  && \Rho \in \H1\calVp \cap \L\infty\calV ,
  \label{regrho}
  \\
  && \Zeta \in \L2\calH .
  \label{regzeta}
\Esist
\Accorpa\Regsoluz regmu regzeta
However, we write $\soluz$ instead of $\jsoluz$ in order to simplify the notation.
As far as the problem under study is concerned, 
we still state it in a weak form as in~\cite{CGS13}, 
by owing to the assumptions \eqref{hpubis} on~$u$.
Namely, we require that
\Bsist
  && \< \dt\Rho , (v,\vG) >_{\calV}
  - \iO \rho u \cdot \nabla v
  + \iO \nabla\mu \cdot \nabla v
  + \iG \nablaG\muG \cdot \nablaG\vG
  = 0
  \non
  \\
  && \quad \hbox{\aet\ and for every $(v,\vG)\in\calV$},
  \label{prima}
  \\
  \separa
  && \tauO \iO \dt\rho \, v
  + \tauG \iG \dt\rhoG \, \vG
  + \iO \nabla\rho \cdot \nabla v
  + \iG \nablaG\rhoG \cdot \nablaG\vG
  \non
  \\
  && \quad {}
  + \iO \bigl( \zeta + \pi(\rho) \bigr) v
  + \iG \bigl( \zetaG + \piG(\rhoG) \bigr) \vG
  = \iO \mu v 
  + \iG \muG \vG
  \non
  \\
  && \quad \hbox{\aet\ and for every $(v,\vG)\in\calV$},
  \label{seconda}
  \\
  && \zeta \in \beta(\rho) \quad \aeQ
  \aand
  \zetaG \in \betaG(\rhoG) \quad \aeS ,
  \label{terza}
  \\
  && \rho(0) = \rhoz
  \quad \aeO \,.
  \label{cauchy}
\Esist
\Accorpa\Pbl prima cauchy

The basic well-posedness result is stated below.
It was proved in \cite[Thm.~2.3]{CGS13}
and holds even in a slightly more general situation.
\juerg{In fact, the condition \eqref{hprhoz} can be weakened, and \eqref{hpbelow} is dispensable}.
Moreover, \eqref{perunicita} is \juerg{needed just for the proof of} uniqueness.

\Bthm
\label{RecallCGS}
Assume that \HPstruttura\ for the structure and \HPdati\ for the data are fulfilled,
and let $T\in(0,+\infty)$ be given.
Then the problem \Pbl\ has a unique solution $\soluz$ satisfying \Regsoluz.
\Ethm

We can obviously draw the following consequence:

\Bcor
\label{Globalsolution}
Assume that \HPstruttura\ and \HPdati\ are satisfied.
Then there exists a unique 6-tuple $\soluz$ defined on $\,(0,+\infty)\,$
that fulfills \Regsoluz\ and solves \Pbl\ for every $T\in(0,+\infty)$.
\Ecor
We remark that the solution is \juerg{actually} much smoother if both $\tauO$ and $\tauG$ are strictly positive
and the initial datum $\rhoz$ satisfies \juerg{further conditions}, as it was proved in \cite[Thm.~2.6]{CGS13}.
\juerg{In particular, under} those assumptions the derivative $\dt\Rho$ can be split \juerg{into components, namely, we
have the representation}
\Beq
  \< \dt\Rho , (v,\vG) >_{\calV}
  = \iO \dt\rho \, v 
  + \iG \dt\rhoG \, \vG
  \label{dtRho}
\Eeq
for every $(v,\vG)\in\calV$, almost everywhere in $(0,+\infty)$.

At this point, given a solution $\soluz$, 
our first aim to investigate its \longtime\ \bhv,
namely, the \omegalimit\
(which we simply term~$\,\omega$, for brevity)
of the component~$\Rho$.
We notice that requiring \eqref{regrho} for every $T\in(0,+\infty)$ 
implies that $\Rho$ is a weakly continuous $\calV$-valued function,
so that the following definition is meaningful.
We~set
\Bsist
  && \omega := \bigl\{
  \Rhoo=\displaystyle\limn\Rho(\tn)
  \quad \hbox{in the weak topology of $\calV$}
  \non
  \\
  && \phantom{\omega := \bigl\{ }
  \mbox{for some sequence $\{\tn\}_{n\in\enne}$ such that $\tn\nearrow+\infty $}
  \bigr\}
  \label{omegalim}
\Esist
and look for the relationship between $\omega$ and the set of stationary solutions
to the system \juerg{which is obtained} from \accorpa{prima}{terza} by ignoring the convective term.
Indeed, assumption \eqref{hpu} implies that $u(t)$ tends to \juerg{the zero function} 
strongly in~$\Lx3$ as $t$ \juerg{approaches} infinity.
It is immediately seen from \eqref{prima} 
that the components $\mu$ and $\muG$
of every stationary solution are \juerg{spatially} constant functions
and that the constant values they assume are the same.
Therefore, \juerg{by} a stationary solution we mean 
a quadruplet $\soluzs$ satisfying, for some $\mus\in\erre$, \juerg{the conditions}
\Bsist
  && \Rhos \in \calV
  \aand
  \Zetas \in \calH,
  \label{regs}
  \\[1mm]
  && \iO \nabla\rhos \cdot \nabla v
  + \iG \nablaG\rhoGs \cdot \nablaG\vG
  + \iO \bigl( \zetas + \pi(\rhos) \bigr) v
  + \iG \bigl( \zetaGs + \piG(\rhoGs) \bigr) \vG
  \non
  \\
  && \quad {}
  = \iO \mus v 
  + \iG \mus \vG
  \quad \hbox{for every $(v,\vG)\in\calV$},
    \label{secondas}
  \\[1mm]
  && \zetas \in \beta(\rhos) \quad \aeO
  \aand
  \zetaGs \in \betaG(\rhoGs) \quad \aeG \,.
  \label{terzas}
\Esist
It is not difficult to show that the conditions \accorpa{regs}{secondas}
imply that the pair $\Rhos$ belongs to $\calW$ and 
satisfies the boundary value problem
\Bsist
  && - \Delta\rhos + \zetas + \pi(\rhos) = \mus
  \quad \aeO ,
  \non
  \\[1mm]
  && \dn\rhos - \DeltaG\rhoGs + \zetaGs + \piG(\rhoGs) = \mus
  \quad \aeG \,.
  \non
\Esist
In \cite{CGS17}, the following result was proved:

\Bthm
\label{Viscouslimit}
Let the assumptions of Corollary~\ref{Globalsolution} be satisfied.
Moreover, assume that the viscosity coefficients $\tauO$ and $\tauG$ are strictly positive,
and let $\soluz$ be the unique global solution on $\,(0,+\infty)$.
Then the \omegalimit\ \eqref{omegalim} is nonempty.
Moreover, for every $\Rhoo\in\omega$,
there exist some $\mus\in\erre$ and a solution $(\rhos,\rhoGs,\zetas,\zetaGs)$ to \accorpa{regs}{terzas}
such that $\Rhoo=\Rhos$.
\Ethm

More precisely, this result holds under slightly weaker assumptions on the velocity field.
The first aim of this paper is to extend Theorem~\ref{Viscouslimit} 
to the pure case, in which $\tauO=\tauG=0$, 
and to the partially viscous situations,
i.e., when either $\tauO>\tauG=0$ or $\tauG>\tauO=0$.
In each of these cases, we need our assumption \eqref{hpu} on~$u$ 
and further conditions on the initial datum.
Thus, we give a list of properties that could be required on~$\rhoz$.
We recall that $\betaeps$ and $\betaGeps$ are the Yosida regularizations
of $\beta$ and~$\betaG$, where $\eta$ is the constant that appears in~\eqref{hpCC}.
Moreover, for $\vG\in\VG$, we use the symbol $\vG^h$ 
for the harmonic extension of $\vG$ to~$\overline\Omega$,~i.e.,
\Beq
  \vG^h \in \Huno , \quad
  \gianni{-\Delta\vG^h=0
  \quad \hbox{in $\Omega$}}
  \aand
  {(\vG^h)}\suG = \vG \,.
  \label{harm}
\Eeq
Here are our possible assumptions:
\Bsist
  && \normaVV{\bigl( -\Delta\rhoz+(\betaeps+\pi)(\rhoz) , \dn\rhoz - \DeltaG\rhoGz+(\betaGeps+\piG)(\rhoGz) \bigr)}\leq C\,,
  \label{hprhozpure}
  \\[1mm]
  && \normaVV{\bigl( (\wG^\eps)^h , \wG^\eps \bigr)}\leq C\,,
  \quad \hbox{where} \quad
  \wG^\eps := \dn\rhoz - \DeltaG\rhoGz+(\betaGeps+\piG)(\rhoGz)\,,
  \qquad
  \label{hprhozOpos}
  \\[1mm]
  && \normaVV{\bigl( -\Delta\rhoz+(\betaeps+\pi)(\rhoz) , (-\Delta\rhoz+(\betaeps+\pi)(\rhoz))\suG \bigr)}\leq C\,,
  \label{hprhozGpos}
\Esist
for some constant $C$ and every $\eps>0$ small enough.
We have chosen to present the above conditions in a unified form,
but some of them could be simplified.
Some comments are given in the forthcoming Remark~\ref{Remhpz}.

\Bthm
\label{Omegalimit}
Suppose that the assumptions of Corollary~\ref{Globalsolution} are satisfied,
and let $(\mu,\muG,\linebreak \rho,\rhoG, \zeta,\zeta_\Gamma)$ be the unique global solution on $(0,+\infty)$.
Moreover, assume either~\eqref{hprhozpure},
or \,$\tauO>0$\, and~\eqref{hprhozOpos},
or \,$\tauG>0$\, and~\eqref{hprhozGpos}.
Then the same \juerg{assertions as in} Theorem~\ref{Viscouslimit}
\juerg{are valid}.
\Ethm

\Brem
\label{Remhpz}
The above \juerg{compatibility} assumptions \accorpa{hprhozpure}{hprhozGpos}
seem to be rather restrictive.
However, one can find \juerg{reasonable} sufficient conditions for them 
\juerg{in the case} that the potentials $f$ and $\fG$ 
are either everywhere defined regular potentials, like \eqref{regpot},
or smooth in $(-1,1)$ and singular at the end-points~$\pm1$, 
like the logarithmic potential~\eqref{logpot}.
Let us consider the latter case.
As for~\eqref{hprhozGpos}, 
by observing that $\normaV{(\vG)^h}\leq\CO\normaVG\vG$ for every $\vG\in\VG$
with $\CO$ depending only on~$\Omega$,
we see that \eqref{hprhozOpos} is equivalent 
to the boundedness in $\VG$ of the second component.
Thus, it is sufficient to assume that 
$\rhoGz\in\HxG2$ (which follows from~\eqref{hprhoz}) and that $\norma\rhoGz_\infty<1$.
Moreover, also the condition \eqref{hprhozGpos} is a consequence of \eqref{hprhoz}
if we assume that $\norma\rhoz_\infty<1$.
So, in fact, \juerg{only the assumption \eqref{hprhozpure} is very restrictive, since 
it postulates that} the second component of the pair appearing there
be the trace of the first one.
This obviously holds if $\pi(0)=\piG(0)=0$ and $\rhoz\in H^3_0(\Omega)$,
since \eqref{hpBeta} implies that $\betaeps(0)=\betaGeps(0)=0$.
If \eqref{hprhoz} is assumed, a~different sufficient condition is the following:
$f$~and $\fG$ are the same logarithmic type potential on~$(-1,1)$,
$\norma\rhoz_\infty<1$, and $(\Delta\rhoz,\DeltaG\rhoGz-\dn\rhoz)\in\calV$.
\Erem

\def\intQ{\int_Q}
\def\intS{\int_\Sigma}
The second aim of this paper is to extend the results of \cite{CGS14}
regarding the control problem to the pure case $\tauO=\tauG=0$
on a fixed time interval~$[0,T]$.
We thus use the simpler notations
\Beq
  Q := \Omega \times (0,T)
  \aand
  \Sigma := \Gamma \times (0,T),
  \non
\Eeq
omitting the subscript~$T$.
The problem addressed in \cite{CGS14} consists of
minimizing the cost functional
\begin{align}
  & \calJ_{gen}(\Mu,\Rho,u)
  := \frac{\beta_1}2 \intQ |\mu-\hmuQ|^2
  + \frac{\beta_2}2 \intS |\muG-\hmuS|^2
  \non
  \\
  & \quad {}
  + \frac{\beta_3}2 \intQ |\rho-\hrhoQ|^2
  + \frac{\beta_4}2 \intS |\rho-\hrhoS|^2
  \non
  \\
  & \quad {}
  + \frac{\beta_5}2 \iO |\rho(T)-\hrhoO|^2
  + \frac{\beta_6}2 \iG |\rhoG(T)-\hrhoG|^2
  + \frac{\beta_7}2 \intQ |u|^2\,,
  \label{gencost} 
\end{align}
subject to the state system \Pbl\ on the time interval $[0,T]$
and to the control constraint $u\in\Uad$, where the convex set $\Uad$ and the related spaces are defined~by
\Bsist
  && \Uad := \bigl\{ u\in\calX: \ |u|\leq\overline U \,\, \aeQ, \ \norma u_{\calX}\leq R_0 \bigr\},
  \label{defUad}
  \\[1mm]
  && \calX := \L2Z \cap (L^\infty(Q))^3 \cap (\H1{\Lx3})^3,
  \label{defcalX}
  \\[1mm]
  && Z := \graffe{ w\in(\Lx2)^3:\ \div w=0 \ \hbox{\,in $\,\Omega$ \ and }\ w\cdot\nu=0 \ \hbox{\,on \,$\Gamma$}}.
  \label{defZ}
\Esist
In \juerg{\eqref{gencost}}, the constants $\beta_i$, $1\le i\le 7$, are nonnegative but not all zero, and $\hmuQ$,
$\hmuS$, $\hrhoQ$, $\hrhoS$, $\hrhoO$, and $\hrhoG$, are prescribed target functions.
In~\eqref{defUad}, the function $\overline U$ and the constant $R_0$ are given
in such a way that $\Uad$ is nonempty.
Moreover, the whole treatement is done in \cite{CGS14} just for  the case $\beta_1=\beta_2=0$.
\juerg{We therefore do the same} in the present paper from the very beginning, for simplicity,
even though such an assumption is not needed immediately.
\juerg{Next, we point out that the results of~\cite{CGS14} were only established 
for potentials $f$ and $\fG$ of logarithmic type}.
Nevertheless, these results also hold true in the case of everywhere defined smooth potentials
whenever a uniform $L^\infty$ bound for the component $\rho$ of the solution can be established.
Indeed, the form of the logarithmic potential is \juerg{(besides its obvious physical importance)}
used there just to ensure that
$\rho$ attains its values in a compact subset of the domain $D(\beta)=D(\betaG)=(-1,1)$
and that the potentials are smooth in such an interval.
Here are our additional assumptions:
\begin{align}
& D(\beta) = D(\betaG) = \erre,
  \aand 
  \hbox{$\beta$, $\betaG$, $\pi$ and $\piG$ are $C^2$ functions}.
  \qquad
  \label{hppotcon}
  \\[1mm]
& \beta_i \geq 0
  \quad \hbox{for $3\leq i\leq 7$}.
  \label{hpcoeff}
  \\[1mm]
& \hrhoQ \in L^2(Q), \quad
  \hrhoS \in L^2(\Sigma) , \quad
  \hrhoO \in \Ldue,
  \aand
  \hrhoG \in \LdueG.
  \label{hptarget}
  \\[1mm]
& \hbox{The function $\overline U\in L^\infty(Q)$ and the constant $R_0>0$ are such that
  } \,\Uad\not=\emptyset.
  \label{hpUad}
\end{align}
We notice that \eqref{hppotcon} allows us to rewrite  \eqref{terza} as
$\zeta=\beta(\rho)$ and $\zetaG=\betaG(\rhoG)$.
Therefore, when speaking of a solution, \juerg{we tacitly assume the validity of these identities and}
just refer to the quadruplet $(\mu,\muG,\rho,\rhoG)$.
We also remark that \eqref{hppotcon} can be expressed in terms of the potentials $f$ and $\fG$ defined in \eqref{defpot}
by just saying that they are $C^3$ functions on the whole real line.
Moreover, since $\beta_1=\beta_2=0$, we rewrite the cost functional \eqref{gencost} in the form
\Bsist
  && \calJ(\Rho,u)
  := \frac{\beta_3}2 \intQ |\rho-\hrhoQ|^2
  + \frac{\beta_4}2 \intS |\rho-\hrhoS|^2
  \non
  \\
  && \phantom{\calJ(\Rho,u) :=}
  + \frac{\beta_5}2 \iO |\rho(T)-\hrhoO|^2
  + \frac{\beta_6}2 \iG |\rhoG(T)-\hrhoG|^2
  + \frac{\beta_7}2 \intQ |u|^2 \,.
  \qquad
  \label{cost} 
\Esist

In the present paper, we can develop a rather complete theory 
only under strong assumptions on the initial datum~$\rhoz$,
for which we require that \eqref{hprhozpure} holds,
even from the very beginning, for simplicity.
We first extend a simplified version of \juerg{the result of} \cite[Thm.~4.1]{CGS14} 
(which required the viscosity coefficients $\tauO$ and~$\tauG$ to be positive)
to~the pure Cahn--Hilliard system.
Indeed, we have the following result.

\Bthm
\label{Optimum}
Assume \HPstruttura\ and \eqref{hppotcon} on the structure,
\HPdati\ and \eqref{hprhozpure} on the data,
\eqref{hpcoeff} and \juerg{\eqref{hptarget}} on the cost functional,
and \eqref{hpUad} on the control box.
Then the problem of minimizing the functional \eqref{cost}
subject to the state system \Pbl\
with $\tauO=\tauG=0$
and to the control constraint $u\in\Uad$
has at least \juerg{one} solution.
\Ethm

The next step is to find necessary conditions for optimality.
We first recall the corresponding results of~\cite{CGS14}.
The main tool is the adjoint problem associated with an optimal control $\ub$
and the corresponding state $(\mub,\muGb,\rhob,\rhoGb)$:
find a quadruplet $(p,\pG,q,\qG)$ satisfying the regularity requirements
\Bsist
  && (p,\pG) \in \L2\calV \,, \quad
  (q,\qG) \in \L\infty\calH \cap \L2\calV\,,
  \label{regpqtau}
  \\
  && (p+\tauO q,\pG+\tauG\qG) \in \H1\calVp\,, 
  \label{regcombinpqtau}
\Esist
and solving 
\Bsist
  && - \< \dt (p+\tauO q,\pG+\tauG\qG) , (v,\vG) >_{\calV}
  + \iO \nabla q \cdot \nabla v
  + \iG \nablaG\qG \cdot \nablaG\vG
  \non
  \\
  && \quad {}
  + \iO \psi q v
  + \iG \psiG \qG \vG
  - \iO \ub \cdot \nabla p \, v
  = \iO \phi_3 \, v
  + \iG \phi_4 \, \vG
  \non
  \\[1mm]
  && \quad \hbox{\aet\ and for every $(v,\vG)\in\calV$},
  \label{primaAvisc}
  \\[2mm]
  \separa
  && \iO \nabla p \cdot \nabla v
  + \iG \nablaG\pG \cdot \nablaG\vG
  = \iO q v 
  + \iG \qG \vG
  \non
  \\[1mm]
  && \quad \hbox{\aet\ and for every $(v,\vG)\in\calV$},
  \label{secondaAvisc}
  \\[2mm]
  && \< (p+\tauO q,\pG+\tauG\qG)(T) , (v,\vG) >_{\calV}
  = \iO \phi_5 \, v
  + \iG \phi_6 \vG
  \non
  \\
  && \quad \hbox{for every $(v,\vG)\in\calV$}, 
  \label{cauchyAvisc}
\Esist
\Accorpa\Adjointtau primaAtau cauchyAtau
where $\tauO$ and $\tauG$ are positive and the following abbreviations are used:
\Bsist
  && \psi := f''(\rhob)
  \aand 
  \psiG := \fG''(\rhoGb) 
  \label{defpsi}
  \\
  && \phi_3 := \beta_3 (\rhob-\hrhoQ) , \quad
  \phi_4 := \beta_4 (\rhoGb-\hrhoS),
  \label{deftrequattro}
  \\
  && \phi_5 := \beta_5 (\rhob(T)-\hrhoO) , \quad
  \phi_6 := \beta_6 (\rhoGb(T)-\hrhoG) .
  \label{defcinquesei}
\Esist
In the quoted paper (see \cite[Thm.~4.4]{CGS14}),
an existence result for the above problem \juerg{was proved
under the assumption that $\tauO$ and $\tauG$ be} positive.
Moreover, the solution \juerg{turned out to be} unique if $\tauO=\tauG$.
In the same paper (see \cite[Thm.~4.6]{CGS14}),
the following optimality condition was derived:
\Beq
  \intQ \bigl( \rhob \, \nabla p + \beta_7 \ub \bigr) \cdot (v-\ub)
  \geq 0
  \quad \hbox{for every $v\in\Uad\,$}.
  \label{necessarytau}
\Eeq
In particular, if $\beta_7>0$, the optimal control
$\ub$ is the $L^2$-projection of $-\frac1{\beta_7}\,\rhob\,\nabla p$ on~$\Uad$.
We once more remark that all this \juerg{was proved in \cite{CGS14} for the case
of logarithmic potentials only; but the same analysis can be carried out for 
everywhere smooth potentials under the assumptions of Theorem~\ref{Optimum}
provided that $\tauO$ and $\tauG$ are positive.} 

In the present paper, we prove a weak version of \juerg{these} results in the pure case,
i.e., \juerg{when} $\tauO=\tauG=0$.
Indeed, \juerg{it turns out that the adjoint problem has to be presented in a time-integrated} form.
To describe it, for every $v\in L^1(Q)$ 
(and similarly for functions in $L^1(\Sigma)$ or in some product space),
we introduce the \juerg{backward-in-time} convolution $1*v$ by the formula
\Beq
  (1*v)(t) := \int_t^T v(s) \, ds
  \quad \aat .
  \label{backconv}
\Eeq
Then, in the pure case,
the adjoint problem associated with an optimal control $\ub$
and the corresponding state $(\mub,\muGb,\rhob,\rhoGb)$ 
consists in finding a quadruplet $(p,\pG,q,\qG)$ 
satisfying the regularity requirement
\Bsist
  && (p,\pG) \in \L2\calV,
  \label{regp}
  \\
  && (q,\qG) \in \L2\calH 
  \quad \hbox{with} \quad
  1*(q,\qG) \in \L\infty\calV,
  \label{regq}
\Esist
and solving the  system
(with the notations \accorpa{defpsi}{defcinquesei})
\Bsist
  && \iO p v
  + \iG \pG \vG
  + \iO \nabla(1*q) \cdot \nabla v
  + \iG \nablaG(1*\qG) \cdot \nabla \vG
  \non
  \\[1mm]
  && \quad {}
  + \iO (1*(\psi q)) \, v
  + \iG (1*(\psiG\qG) \, \vG 
  + \iO (1*(\ub\cdot\nabla p)) \, v
  \non
  \\[1mm]
  && {} = \iO (1*\phi_3) v
  + \iG (1*\phi_4) \vG
  + \iO \phi_5 v
  + \iG \phi_6  \vG
  \non
  \\[2mm]
  && \quad \hbox{\aat\ and for every $(v,\vG)\in\calV$}
  \label{primaA}
  \\[2mm]
  && \intQ \nabla p \cdot \nabla v
  + \intS \nablaG\pG \cdot \nablaG\vG
  = \intQ q v 
  + \intS \qG \vG
  \non
  \\[2mm]
  && \quad \hbox{\aat\ and for every $(v,\vG)\in\L2\calV$}.
  \label{secondaA}
\Esist

Here is our result on \juerg{first-order} optimality conditions.

\Bthm
\label{Necessary}
Suppose that the conditions \HPstruttura\ and \eqref{hppotcon} on the structure,
\HPdati\ and \eqref{hprhozpure} on the data,
\eqref{hpcoeff} and \eqref{hptarget} on the cost functional,
and \eqref{hpUad} on the control box, are fulfilled.
Moreover, assume that $\,\ub\,$ and $(\mub,\muGb,\rhob,\rhoGb)$ are
an optimal control and the corresponding optimal state, respectively.
Then there exists at least one quadruplet $(p,\pG,q,\qG)$
satisfying the conditions \accorpa{regp}{regq} and solving both \accorpa{primaA}{secondaA}
and the variational inequality
\Beq
  \intQ (\rhob \, \nabla p + \beta_7 \ub) \cdot (v-\ub) \geq 0
  \quad \hbox{for every $v\in\Uad$}.
  \label{necessary}
\Eeq
In particular, if $\,\beta_7>0$, then the optimal control
\,$\ub$\, is the $L^2$-projection of $\,-\frac1{\beta_7}\,\rhob\,\nabla p\,$ \gianni{on~$\,\Uad\,$.}
\Ethm

\Brem
\label{UniqA}
We cannot prove a uniqueness result for the solution to problem~\accorpa{primaA}{secondaA}, unfortunately.
However, as will be stated in the forthcoming Remark~\ref{RemuniqA},
the solution is unique provided it is somewhat smoother.
Namely, it is needed that its component $(q,\qG)$ belongs to~$\L2\calV$,
\juerg{which  we are not able to prove}
(while the regularity $(p,\pG)\in\L2\calW$ we did not require
is true and immediately follows from \eqref{secondaA} and the first \eqref{regq}
by applying \cite[Lem.~3.1]{CGS13}).
\Erem

As \juerg{it has been stated} at the end of the introduction,
the proofs of our results will be given in Sections~\ref{LONGTIME} and~\ref{CONTROL}.
\juerg{In order to be able to carry out these proofs, we will now introduce some auxiliary tools, namely, 
the generalized mean value, the related spaces, 
and the operator~$\calN$. In doing this, we will be very brief, referring}  
to \cite[Sect.~2]{CGS13} for further details.
We set
\Beq
  \mean\gstar := \frac { \< \gstar,(1,1) >_{\calV} } { |\Omega| + |\Gamma| }
  \quad \hbox{for $\gstar\in\calVp$}
  \label{genmean}
\Eeq
and observe~that
\Beq
  \mean(v,\vG) = \frac { \iO v + \iG \vG } { |\Omega| + |\Gamma| }
  \quad \hbox{if $(v,\vG)\in\calH$} \,.
  \label{usualmean}
\Eeq
Notice that the constant $\mz$ \juerg{appearing in} assumption \eqref{hpmz} 
is nothing but the mean value \,$\mean(\rhoz,\rhoGz)$,
and that taking $(v,\vG)=(|\Omega|+|\Gamma|)^{-1}(1,1)$ in \eqref{seconda}
yields the conservation property for the component $\Rho$ of the solution,
\Beq
  \dt \mean\Rho = 0,
  \quad \hbox{whence} \quad
  \mean \Rho(t) = \mz
  \quad \hbox{for every $t\in[0,T]$}.
  \label{conservation}
\Eeq
We also stress that the function
\Beq
  \calV \ni (v,\vG) \mapsto
  \bigl( \norma{\nabla v}_{\L2H}^2 + \norma{\nablaG\vG}_{\L2\HG}^2 + |\mean(v,\vG)|^2 \bigr)^{1/2}
  \label{normacalV}
\Eeq
yields a Hilbert norm on $\calV$ that is equivalent to the natural one.
Now, we set
\Beq
  \calVsz := \graffe{ \gstar\in\calVp : \ \mean\gstar = 0 } , \quad
  \calHz := \calH \cap \calVsz
  \aand
  \calVz := \calV \cap \calVsz, 
  \label{defcalVz}
\Eeq
\Accorpa\Defspazi defspaziO defcalVz
and notice that the function
\Beq
  \calVz \ni (v,\vG) \mapsto \norma{(v,\vG)}_{\calVz}
  := \bigl( \norma{\nabla v}_{\L2H}^2 + \norma{\nablaG\vG}_{\L2\HG}^2 \bigr)^{1/2}
  \label{normaVz}
\Eeq
is a Hilbert norm on $\calVz$ which is equivalent to the usual one.
Next, we define the operator $\calN:\calVsz\to\calVz$ 
(which~will be applied to $\calVsz$-valued functions as well)
as~follows: for every element $\gstar\in\calVsz$, we have that
\Bsist
  && \hbox{$\calN\gstar=(\calNO\gstar,\calNG\gstar)$ is the unique pair $\Xi\in\calVz$ such that}
  \non
  \\
  \noalign{\smallskip}
  && \iO \nabla\xi \cdot \nabla v + \iG \nablaG\xiG \cdot \nablaG\vG
  = \< \gstar , (v,\vG) >_{\calV}
  \quad \hbox{for every $(v,\vG)\in\calV$}.
  \qquad
  \label{defN}
\Esist
It turns out that $\calN$ is well defined, linear, symmetric, and bijective.
Therefore, if we set
\Beq
  \norma\gstar_* := \norma{\calN\gstar}_{\calVz},
  \quad \hbox{for $\gstar\in\calVsz$},
  \label{normastar}
\Eeq
then we obtain a Hilbert norm on $\calVsz$
(equivalent to the norm induced by the norm of~$\calVp$),
and we have that
\Beq
  \< \gstar , \calN\gstar >_{\calV}
  = \norma\gstar_*^2
  \quad \hbox{for every $\gstar\in\calVsz$}.
  \label{propN}
\Eeq
Furthermore, we notice that
\Beq
  \< \dt\gstar , \calN\gstar >_{\calV}
  =  \frac 12 \, \frac d{dt} \, \norma\gstar_*^2
  \quad \aet,
  \quad \hbox{for every $\gstar \in \H1\calVsz$}.
  \label{propNdta}
\Eeq
It is easy to see that $\calN\gstar$ belongs to $\calW$ whenever $\gstar\in\calHz$ 
and that
\Beq
  \normaWW{\calN\gstar} \leq \CO \normaHH\gstar
  \quad \hbox{for every $\gstar\in\calHz$},
  \label{regN}
\Eeq
where $\CO$ depends only on $\Omega$.

Besides the above tools, we will repeatedly use the Young inequality
\Beq
  a\,b \leq \delta\,a^2 + \frac 1{4\delta} \, b^2
  \quad \hbox{for all $a,b\in\erre$ and $\delta>0$},
  \label{young}
\Eeq
as well as H\"older's inequality and the Sobolev inequality
\Beq
  \norma v_p
  \leq \CO \normaV v
  \quad \hbox{for every $p\in[1,6]$ and $v\in V$},
  \label{sobolev}
\Eeq
which is related to the continuous embedding $V\subset L^p(\Omega)$
for $p\in[1,6]$
(since $\Omega$ is three-dimensional, bounded and smooth).
In~particular, by also using the equivalent norm \eqref{normacalV} on~$\calV$,
we have that
\Beq
  \norma v_6^2
  \leq \CO \bigl( \norma{\nabla v}_{\L2H}^2 + \norma{\nablaG\vG}_{\L2\HG}^2 + |\mean(v,\vG)|^2 \bigr)
  \label{sobolevbis}
\Eeq
for every $(v,\vG)\in\calV$.
In both \eqref{sobolev} and \eqref{sobolevbis}, the constant $\CO$ depends only on~$\Omega$.
We also account for the compact embedding $\calV\subset\calH$
and for the corresponding compactness inequality
\Beq
  \normaHH{(v,\vG)}^2
  \leq \delta \, \normaHH{(\nabla v,\nablaG\vG)}^2
  + C_{\delta,\Omega} \, \norma{(v,\vG)}_{\calVp}^2
  \quad \hbox{for every $(v,\vG)\in\calV$},
  \label{compact}
\Eeq
where $\delta>0$ is arbitrary,
and where the constant $C_{\delta,\Omega}$ depends only on $\Omega$ and~$\delta$.

Finally, as far as constants are concerned,
we \juerg{employ the following general rule: the} small-case symbol $c$ stands for different 
constants which depend only
on~$\Omega$, the structure of our system 
and the norms of the data involved in the assumptions \HPdati.
A~notation like~$c_\delta$ (in particular, with $\delta=T$) 
allows the constant to depend on the positive
parameter~$\delta$, in addition.
Hence, the meaning of $c$ and $c_\delta$ might
change from line to line and even within the same chain of inequalities.  
On the contrary, we mark those constants that we want to refer~to
by using a different notation (e.g., a~capital letter).


\section{\Longtime\ \bhv}
\label{LONGTIME}
\setcounter{equation}{0}

This section is devoted to the proof of Theorem~\ref{Omegalimit}.
However, the procedure we follow is useful for the next section as well, 
where we prove our results concerning the control problem in the pure case.
Here, we fix any global solution $\soluz$ once and for all.
Our method closely follows the proof of Theorem~\ref{Viscouslimit} performed in~\cite{CGS17};
it thus relies on some global a~priori estimates
and on the study of the \bhv\ of the a solution on intervals 
of a fixed length~$T$ whose endpoints \juerg{approach} infinity.
However, the proof of each crucial estimate \juerg{will have} to be modified.

For brevity, we often proceed formally.
In particular, we behave as if the solution were smooth
and use the identity
\Beq
  \< \dt\Rho , (v,\vG) >_{\calV}
  = \iO \dt\rho \, v 
  + \iG \dt\rhoG \, \vG  
  \quad \hbox{for every $(v,\vG)\in\calV$},
  \label{dtRhobis}
\Eeq
which is not justified, since \juerg{this representation holds true} only in the viscous case $\tauO>0$, $\tauG>0$ (see \eqref{dtRho}). Moreover, we \juerg{argue as if} even the graphs $\beta$ and $\betaG$ were everywhere defined 
smooth functions
and often write $\beta(\rho)$ and $\betaG(\rhoG)$ instead of $\zeta$ and~$\zetaG$, respectively.
However, the estimates obtained in this way \juerg{can} be performed rigorously on a proper approximating problem,
uniformly with respect to the approximation parameter.
The best choice \juerg{for such an approximating problem could} be one of the following:
$a)$ \,\,the \,$\eps$-problem \juerg{which is} analogous to the one introduced in \cite{CGS13} 
and obtained by replacing the graphs $\beta$ and $\betaG$ 
with very smooth functions $\betaeps$ and~$\betaGeps$
that behave like the Yosida regularizations
(i.e., with similar boundedness and convergence properties,
like the $C^\infty$ approximations introduced in \cite[Sect.~3]{GiRo});
$b)$ \,\,the Faedo--Galerkin scheme 
(depending on the parameter $n\in\enne$) 
used in \cite{CGS13} \juerg{in order} to discretize and then solve the $\eps$-problem:
indeed, all of the components of its solution are very smooth.

Now, we start proving some global estimates,
assuming that $\tauO$ and $\tauG$ are nonnegative
in order to cover all relevant situations at the same time, if possible.
However, we will sometimes be forced to distinguish 
between the pure case and the partially viscous ones.

\step
First global estimate

We recall the conservation property \eqref{conservation},
write the equations \eqref{prima} and \eqref{seconda} at the time~$s$,
and test them by $\Mu(s)+\calN(\dt\Rho(s))\in\calV$ and $2\dt\Rho(s)\in\calV$, respectively.
Then, we use \eqref{propN} with $\gstar=\dt\Rho(s)$,
integrate with respect to $s$ over $(0,t)$ with an arbitrary~$t>0$, 
sum up and rearrange. We obtain the identity
\Bsist
  && \intQt \dt\rho \, \mu
  + \intSt \dt\rhoG \, \muG
  + \intQt |\nabla\mu|^2
  + \intSt |\nablaG\muG|^2  
  \non
  \\
  && \quad {}
  + \iot \norma{\dt\Rho(s)}_*^2 \, ds
  + \intQt \nabla\mu \cdot \nabla\calNO(\dt\Rho)
  + \intSt \nablaG\muG \cdot \nablaG\calNG(\dt\Rho)
  \non
  \\
  \separa
  && \quad {}
  + 2\tauO \intQt |\dt\rho|^2
  + 2\tauG \intSt |\dt\rhoG|^2  
  \non
  \\
  && \quad {}
  + \iO |\nabla\rho(t)|^2
  + \iG |\nablaG\rhoG(t)|^2
  + 2 \iO f(\rho(t))
  + 2 \iG \fG(\rhoG(t))
  \non
  \\
  \separa
  && = \intQt \rho u \cdot \bigl( \nabla\mu + \nabla\calNO(\dt\Rho) \bigr)
  + 2 \intQt \mu \dt\rho
  + 2 \intSt \muG \dt\rhoG
  \non
  \\
  && \quad {}
  + \iO |\nabla\rhoz|^2
  + \iG |\nablaG\rhoGz|^2
  + 2 \iO f(\rhoz)
  + 2 \iG \fG(\rhoGz).
  \non
\Esist
Some integrals cancel out by the definition \eqref{defN} of~$\calN$,
the terms on the \lhs\ containing $f$ and $\fG$ are bounded from below by~\eqref{hpbelow},
and the ones on the \rhs\ involving the initial values are finite by~\eqref{hprhoz}.
We deal with the convective term 
by owing to the Young and \Holder\ inequalities,
the Sobolev type inequality \eqref{sobolevbis}, and the conservation property~\eqref{conservation}.
We~have
\Bsist
  && \intQt \rho u \cdot \bigl( \nabla\mu + \nabla\calNO(\dt\Rho) \bigr)
  \non
  \\
  && \leq \frac 12 \intQt |\nabla\mu|^2
  + \frac 12 \iot \norma{\calN(\dt\Rho(s))}_{\calVz}^2 \, ds
  + \iot \norma{u(s)}_3^2 \, \norma{\rho(s)}_6^2 \, ds
  \non 
  \\
  && \leq \frac 12 \intQt |\nabla\mu|^2
  + \frac 12 \iot \norma{\dt\Rho(s)}_*^2 \, ds
  \non
  \\
  && \quad {}
  + c \iot \norma{u(s)}_3^2 \, \bigl( \norma{\nabla\rho(s)}_{\L2H}^2 + \norma{\nablaG\rhoG(s)}_{\L2\HG}^2 + \mz^2 \bigr) \, ds .
  \non
\Esist
Since the function $s\mapsto\norma{u(s)}_3^2$ belongs to $L^1(0,+\infty)$ by~\eqref{hpu},
we can apply the Gronwall lemma on $(0,+\infty)$ and obtain~that
\Bsist
  && \intQi |\nabla\mu|^2
  + \intSi |\nablaG\muG|^2
  + \int_0^{+\infty} \norma{\dt\Rho(t)}_*^2\,dt
  < + \infty\,,
  \label{globint}
  \\
  \separa
  && \intQi |\dt\rho|^2 < +\infty
  \quad \hbox{if $\tauO>0$}
  \aand
  \intSi |\dt\rhoG|^2 < +\infty
  \quad \hbox{if $\tauG>0$}\,,
  \label{globinttau}
  \\
  && f(\rho) \in \LL\infty\Luno
  \aand
  \fG(\rhoG) \in \LL\infty\LunoG\,,
  \label{stimef}
\Esist
as well as $(\nabla\rho,\nablaG\rhoG)\in(\LL\infty\calH)^3$.
From this, by accounting for the conservation property~\eqref{conservation} once more,
we conclude that
\Beq
  \Rho \in \LL\infty\calV .
  \label{primastima}
\Eeq

\step
Consequence

By using the quadratic growth of $\Pi$ and $\PiG$, \juerg{which is} 
implied by the \Lip\ continuity of their derivatives,
and combining with \eqref{stimef} with~\eqref{primastima},
we deduce that
\Beq
  \Beta(\rho) \in  \LL\infty\Luno
  \aand
  \BetaG(\rhoG) \in \LL\infty\LunoG \,.
  \label{stimeBeta}
\Eeq

\step
Second global estimate

We formally differentiate the equations \eqref{prima} and \eqref{seconda}
with respect to time, behaving as if $\beta$ and $\betaG$ were smooth functions
and writing $\beta(\rho)$ and $\betaG(\rhoG)$
instead of $\zeta$ and~$\zetaG$ (see~\eqref{terza}).
We obtain 
\Bsist
  && \iO \dt^2\rho \, v
  + \iG \dt^2\rhoG \, \vG
  + \iO \nabla\dt\mu \cdot \nabla v
  + \iG \nablaG\dt\muG \cdot \nablaG\vG
  \non
  \\
  && = \iO \bigl( \dt\rho \, u + \rho \dt u \bigr) \cdot \nabla v\,,
  \label{dtprima}
  \\
  \separa
  && \tauO \iO \dt^2\rho \, v
  + \tauG \iG \dt^2\rhoG \, \vG 
  + \iO \nabla\dt\rho \cdot \nabla v
  + \iG \nablaG\dt\rhoG \cdot \nablaG\vG
  \non
  \\
  && \quad {}
  + \iO \beta'(\rho) \dt\rho \, v
  + \iG \betaG'(\rhoG) \dt\rhoG \, \vG
  \non
  \\
  && = \iO \dt\mu \, v 
  + \iG \dt\muG \, \vG
  - \iO \pi'(\rho) \dt\rho \, v
  - \iG \piG'(\rhoG) \dt\rhoG \, \vG \,,
  \label{dtseconda}
\Esist
\Aet,  and for every $(v,\vG)\in\calV$.
Recalling that $\dt\Rho$ is $\calVz$-valued by \eqref{conservation},
so that $\calN\dt\Rho$ is well defined,
we write the above equations at the time~$s$
and test them by $\calN\dt\Rho(s)$ and $\dt\Rho(s)$, respectively.
Then we integrate with respect to $s$ over $(0,t)$ with an arbitrary $t>0$ and sum~up.
It follows that
\Bsist
  && \iot \< \dt^2\Rho(s) , \calN\dt\Rho(s) >_\calV \, ds
  \non
  \\
  && \quad {}
  + \intQt \nabla\dt\mu \cdot \nabla\calNO(\dt\Rho)
  + \intSt \nablaG\dt\muG \cdot \nablaG\calNG(\dt\Rho)
  \non
  \\
  && \quad {}
  + \tauO \intQt \dt^2\rho \, \dt\rho
  + \tauG \intSt \dt^2\rhoG \, \dt\rhoG
  + \intQt |\nabla\dt\rho|^2
  + \intSt |\nablaG\dt\rhoG|^2
  \non
  \\
  && \quad {}
  + \intQt \beta'(\rho) |\dt\rho|^2
  + \intSt \betaG'(\rhoG) |\dt\rhoG|^2
  \non
  \\
  \separa
  && = \intQt \bigl( \dt\rho \, u + \rho \dt u \bigr) \cdot \nabla\calNO(\dt\Rho)
  \non
  \\
  && \quad {}
  + \intQt \dt\mu \, \dt\rho 
  + \intSt \dt\muG \, \dt\juerg{\rhoG}
  - \intQt \pi'(\rho) |\dt\rho|^2
  - \intSt \piG'(\rhoG) |\dt\rhoG|^2 .
  \non
\Esist
The integrals containing $\dt\mu$ and $\dt\muG$ cancel out
by the definition \eqref{defN} of~$\calN$
(with the choices $\gstar=\dt\Rho(s)$ and $(v,\vG)=\dt\Mu(s)$),
and the terms involving $\beta'$ and $\betaG'$ are nonnegative.
Moreover, the last two integrals on the \rhs\ can be dealt with
by accounting for the \Lip\ continuity of $\pi$ and~$\piG$
and the compactness inequality~\eqref{compact},
\Bsist
  && - \intQt \pi'(\rho) |\dt\rho|^2
  - \intSt \piG'(\rhoG) |\dt\rhoG|^2
  \non
  \\
  && \leq \frac 12 \intQt |\nabla\dt\rho|^2
  + \frac 12 \intSt |\nablaG\dt\rhoG|^2
  + c \iot \norma{\dt\Rho(s)}_*^2 \, ds \,.
  \non
\Esist
Since the last integral is bounded by~\eqref{globint},
coming back to the previous identity and owing to \eqref{propNdta} for the first term on the \lhs,
we deduce that
\Bsist
  && \frac 12 \, \norma{\dt\Rho(t)}_*^2
  + \frac \tauO 2 \iO |\dt\rho(t)|^2
  + \frac \tauG 2 \iG |\dt\rhoG(t)|^2
  + \frac 12 \intQt |\nabla\dt\rho|^2
  + \frac 12 \intSt |\nablaG\dt\rhoG|^2
  \non
  \\
  && \leq \frac 12 \, \norma{\dt\Rho(0)}_*^2
  + \frac \tauO 2 \iO |\dt\rho(0)|^2
  + \frac \tauG 2 \iG |\dt\rhoG(0)|^2
  \non
  \\
  && \quad {}
  + \intQt \bigl( \dt\rho \, u + \rho \dt u \bigr) \cdot \nabla\calNO(\dt\Rho)
  + c \,.
  \label{persecondastima}
\Esist
Thus, it suffices to estimate the last integral and obtain bounds for 
the time derivatives evaluated at~$0$.
For the first aim, we use the \Holder, Sobolev and Young inequalities (in particular, \eqref{sobolevbis}),
the conservation property~\eqref{conservation},
and the already established estimates \eqref{globint} and \eqref{primastima}.
We then have
\Bsist
  && \intQt \bigl( \dt\rho \, u + \rho \dt u \bigr) \cdot \nabla\calNO(\dt\Rho)
  \non
  \\
  && \leq \iot \bigl(
    \norma{\dt\rho(s)}_6 \, \norma{u(s)}_3  
    + \norma{\rho(s)}_6 \, \norma{\dt u(s)}_3  
  \bigr) \norma{\nabla\calNO(\dt\Rho(s))}_2 \, ds
  \non
  \\
  && \leq \frac 14 \, \Bigl( \intQt |\nabla\dt\rho|^2 + \intSt |\nablaG\dt\rhoG|^2 \Bigr)
  \non
  \\
  && \quad {}
  + c \, \norma u_{\LL\infty{\Lx3}}^2 \iot \norma{\dt\Rho(s)}_*^2 \, ds
  + c \, \norma\rho_{\LL\infty V}^2 \iot \norma{\dt u(s)}_3^2 \, ds  
  \non
  \\
  && \leq \frac 14 \, \Bigl( \intQt |\nabla\dt\rho|^2 + \intSt |\nablaG\dt\rhoG|^2 \Bigr)
  + c \,.
  \non
\Esist
In order to find bounds for the initial value of the time derivatives,
we have to distinguish between the pure case and the partially \juerg{viscous} ones.
Assume first that $\tauO=\tauG=0$.
Then the only initial value that appears on the \rhs\ of \eqref{persecondastima}
is $\norma{\dt\Rho(0)}_*$,
and we can estimate it by accounting for \eqref{hprhozpure},
which we formally write~as
\Beq
  \bigl( -\Delta\rhoz+(\beta+\pi)(\rhoz) , \dn\rhoz - \DeltaG\rhoGz+(\betaG+\piG)(\rhoGz) \bigr)  
  \in \calV \,. 
  \non
\Eeq
We write \eqref{prima} and \eqref{seconda} at the time $t=0$,
test them by $\Mu(0)+\calN(\dt\Rho(0))$ and $2\dt\Rho(0)$, respectively, and sum~up.
Then, we owe to \eqref{propN} with $\gstar=\dt\Rho(0)$,
integrate by parts the terms involving $\nabla\rhoz$ and~$\nablaG\rhoGz$,
and rearrange.
We obtain
\Bsist
  && \iO \dt\rho(0) \, \mu(0) 
  + \iG \dt\rhoG(0) \, \muG(0)
  + \iO |\nabla\mu(0)|^2
  + \iG |\nablaG\muG(0)|^2
  \non
  \\
  && \quad {}
  + \norma{\dt\Rho(0)}_*^2
  \non
  \\
  && \quad {}
  + \iO \nabla\mu(0) \cdot \nabla\calNO(\dt\Rho(0))
  + \iG \nablaG\muG(0) \cdot \nablaG\calNG(\dt\Rho(0))
  \non
  \\
  \separa
  && = \iO \rhoz u(0) \cdot \bigl( \nabla\mu(0) + \nabla\calNO(\dt\Rho(0)) \bigr)
  + 2 \iO \mu(0) \dt\rho(0)
  + 2 \iG \muG(0) \dt\rhoG(0)
  \non
  \\
  && \quad {}
  - 2 \iO \bigl( -\Delta\rhoz+(\beta+\pi)(\rhoz) \bigr) \, \dt\rho(0)
  \non
  \\
  && \quad {}
  - 2 \iG \bigl( \dn\rhoz - \DeltaG\rhoGz+(\betaG+\pi)(\rhoGz) \bigr) \, \dt\rhoG(0).
  \qquad\enskip
  \label{perstimadtiniz}
\Esist
Some integrals cancel each other, either trivially or by the definition of~$\calN$.
In particular, what remains on the \rhs\ are just the first and the last two terms,
the sum of which is estimated from above~by
\Beq
  \frac 12 \, \norma{\dt\Rho(0)}_*^2
  + c \normaVV{\bigl( -\Delta\rhoz+(\beta+\pi)(\rhoz) , \dn\rhoz - \DeltaG\rhoGz+(\betaG+\piG)(\rhoGz) \bigr)}^2, 
  \non
\Eeq
thanks to the Young inequality.
Thus, we infer~that
\Beq
  \iO |\nabla\mu(0)|^2
  + \frac 12 \, \norma{\dt\Rho(0)}_*^2
  \leq \iO \rhoz u(0) \cdot \bigl( \nabla\mu(0) + \nabla\calNO(\dt\Rho(0)) \bigr)
  + c \,.
  \non
\Eeq
On the other hand, by the \Holder\ inequality 
and due to our assumptions \eqref{hprhoz} and \eqref{hpu} on $\rhoz$ and~$u$, we also have that
\Bsist
  && \iO \rhoz u(0) \cdot \bigl( \nabla\mu(0) + \nabla\calNO(\dt\Rho(0)) \bigr)
  \non
  \\
  && \leq \frac 12 \iO |\nabla\mu(0)|^2
  + \frac 14 \, \norma{\nabla\calNO(\dt\Rho(0))}_2^2
  + c \, \bigl( \norma\rhoz_6^2 + \norma{u(0)}_3^2 \bigr)
  \non
  \\
  && \leq \frac 12 \iO |\nabla\mu(0)|^2
  + \frac 14 \, \norma{\calN(\dt\Rho(0))}_{\calVz}^2
  + c 
  \non
  \\
  && \leq \frac 12 \iO |\nabla\mu(0)|^2
  + \frac 14 \, \norma{\dt\Rho(0)}_*^2
  + c \,.
  \non
\Esist
By combining this \gianni{estimate} with the above inequality, 
we obtain the \juerg{sought} bound for~$\norma{\dt\Rho(0)}_*$.

If one of the viscosity parameters $\tauO$ and $\tauG$ is positive,
a~similar argument applies.
Moreover, we can owe to the weaker assumption \eqref{hprhoz}
and to either \eqref{hprhozOpos} or \eqref{hprhozGpos}
(instead of \eqref{hprhozpure}) on~the initial datum~$\rhoz$,
as we show at once.
The \rhs\ of \eqref{persecondastima} also contains the integrals
\Beq
  \frac \tauO 2 \iO |\dt\rho(0)|^2
  + \frac \tauG 2 \iO |\dt\rhoG(0)|^2\,,
  \label{dtrhozvisc}
\Eeq
and one of them is significant and has to be estimated as well.
Assume first \juerg{that} \,$\tauO>0$ and $\tauG=0$.
Then, testing \eqref{seconda} for\, $t=0$\, by \,$\dt\Rho(0)$, as \juerg{done above},
yields the additional term $\tauO\iO|\dt\rho(0)|^2$
on the \lhs\ of~\eqref{perstimadtiniz}, whence \eqref{dtrhozvisc} can actually be estimated.
But the presence of $\tauO\iO|\dt\rho(0)|^2$ on the \lhs\
also helps in estimating the last two terms on the \rhs\ of \eqref{perstimadtiniz}
by taking advantage of \eqref{hprhoz} and \eqref{hprhozOpos} instead of the more restrictive~\eqref{hprhozpure}.
We formally write \eqref{hprhozOpos}~as
\Beq
  (\wG^h,\wG) \in \calV
  \quad \hbox{where} \quad
  \wG := \dn\rhoz-\DeltaG\rhoGz+(\betaG+\piG)(\rhoGz) .
  \non
\Eeq
Since $\normaH{\vG^h}\leq c\normaVV{(\vG^h,\vG)}$ for every $\vG\in\VG$,
we have 
\Bsist
  && - 2 \iO \bigl( -\Delta\rhoz+(\beta+\pi)(\rhoz) \bigr) \, \dt\rho(0)
  - 2 \iG \bigl( \dn\rhoz-\DeltaG\rhoGz+(\betaG+\piG)(\rhoGz) \bigr) \, \dt\rhoG(0)
  \non
  \\
  && = - 2 \iO \bigl[
    \bigl( -\Delta\rhoz+(\beta+\pi)(\rhoz) \bigr)
    - \wG^h
  \bigr] \, \dt\rho(0)
  - 2 \, \< \dt\Rho(0) , (\wG^h,\wG) >_{\calV}
  \non
  \\
  && \leq \frac \tauO 4 \iO |\dt\rho(0)|^2
  + \frac 14 \, \norma{\dt\Rho(0)}_*^2
  \non
  \\
  && \quad {}
  + c \, \normaH{-\Delta\rhoz+(\beta+\pi)(\rhoz)}^2 
  + c \, \normaH{\wG^h}^2
  + c \, \normaVV{(\wG^h , \wG)}
  \non
  \\
  && \leq \frac \tauO 4 \iO |\dt\rho(0)|^2
  + \frac 14 \, \norma{\dt\Rho(0)}_*^2
  + c \,,
  \non
\Esist
and this leads to the desired estimate.

Assume now \juerg{that} \,$\tauG>0$ and $\tauO=0$.
Similarly as before, 
testing \eqref{dtseconda} by $\dt\Rho$ yields a nonnegative contribution to the \lhs\
and the second term of \eqref{dtrhozvisc} on the \rhs.
But testing \eqref{seconda} for $t=0$ by $\dt\Rho(0)$
now produces the integral $\tauG\iG|\dt\rhoG(0)|^2$
on the \lhs\ of \eqref{perstimadtiniz}
(whence \eqref{dtrhozvisc} can be estimated)
and allows us to treat the last two terms of the \rhs\ of \eqref{perstimadtiniz}
by owing to \eqref{hprhoz} and \eqref{hprhozGpos}.
We formally write the latter~as
\Beq
  \bigl( -\Delta\rhoz+(\beta+\pi)(\rhoz) , (-\Delta\rhoz+(\beta+\pi)(\rhoz))\suG \bigr) \in \calV \,.
  \non
\Eeq
By taking advantage of the inequalities $\normaHG\vG\leq c\,\normaVV{(v,\vG)}$ for every $(v,\vG)\in\calV$, 
and $\normaHG{\dn v}\leq c\,\normaW v$ for every $v\in W$
(the former is obvious and the latter is well known), 
we have the estimate
\Bsist
  && - 2 \iO \bigl( -\Delta\rhoz+(\beta+\pi)(\rhoz) \bigr) \, \dt\rho(0)
  - 2 \iG \bigl( \dn\rhoz-\DeltaG\rhoGz+(\betaG+\pi)(\rhoGz) \bigr) \, \dt\rhoG(0)
  \non
  \\
  && = -2 \, \< \dt\Rho(0) , \bigl( -\Delta\rhoz+(\beta+\pi)(\rhoz) , (-\Delta\rhoz+(\beta+\pi)(\rhoz))\suG \bigr) >_{\calV}
  \non
  \\
  && \quad {}
  + 2 \iG \bigl[
    (-\Delta\rhoz+(\beta+\pi)(\rhoz))\suG
    - (-\DeltaG\rhoGz+(\betaG+\piG)(\rhoGz))
    - \dn\rhoz
  \bigr] \dt\rhoG(0)
  \non
  \\
  && \leq \frac 14 \, \norma{\dt\Rho(0)}_*^2
  + \frac \tauG 4 \iG |\dt\rhoG(0)|^2
  \non  
  \\
  && \quad {}
  + c \, \normaVV{\bigl( -\Delta\rhoz+(\beta+\pi)(\rhoz) , (-\Delta\rhoz+(\beta+\pi)(\rhoz))\suG \bigr)}^2
  \non
  \\
  && \quad {}
  + c \, \normaHG{(-\Delta\rhoz+(\beta+\pi)(\rhoz))\suG}^2
  + c \, \normaHG{-\DeltaG\rhoGz+(\betaG+\piG)(\rhoGz)}^2
  + c \, \normaHG{\dn\rhoz}^2
  \non
  \\
  && \leq \frac 14 \, \norma{\dt\Rho(0)}_*^2
  + \frac \tauG 4 \iG |\dt\rhoG(0)|^2
  + c\,,
  \non  
\Esist
whence the desired bounds for the initial values of the time derivatives follow.

Now, coming back to \eqref{persecondastima}, 
and taking these estimates into account, we conclude that
\Bsist
  && \dt\Rho \in \LL\infty\calVp \,,
  \label{secondastima}
  \\
  && \dt\rho \in \LL\infty H
  \quad \hbox{if $\tauO>0$}
  \aand
  \dt\rhoG \in \LL\infty\HG
  \quad \hbox{if $\tauG>0$}\,.
  \qquad
  \label{secondastimatau}
\Esist
As a by-product, we also obtain the less important result that \,$\dt\Rho\in\LL2\calV$.

\medskip

In the rest of the section, we only consider the pure case $\tauO=\tauG=0$.
However, the whole argument works in the partially viscous cases as well.
Indeed, testing equations as we do 
would produce just additional contributions that can be dealt with in a trivial way.

\step
Third global estimate

Once again, we write $\zeta=\beta(\rho)$ and $\zetaG=\beta(\rhoG)$
and first notice that the inclusion $D(\betaG)\subseteq D(\beta)$ (see~\eqref{hpCC}) 
and assumption \eqref{hpmz} imply that
\Beq
  \beta(r) (r-\mz)
  \geq \delta_0 |\beta(r)| - C_0
  \aand
  \betaG(r) (r-\mz)
  \geq \delta_0 |\betaG(r)| - C_0
  \label{trickMZ}
\Eeq
for every $r$ belonging to the respective domains,
where $\delta_0$ and $C_0$ are some positive constants that depend only on $\beta$, $\betaG$
and on the position of $\mz$ in the interior of~$D(\betaG)$ and of~$D(\beta)$
(see, e.g. \cite[p.~908]{GiMiSchi}).
Now, we recall the conservation property \eqref{conservation}
and test \eqref{prima} and \eqref{seconda} by 
$\calN(\rho-\mz,\rhoG-\mz)$ and $(\rho-\mz,\rhoG-\mz)$, respectively.
Then we add, without integrating with respect to time.
We obtain, \Aet,
\Bsist
  && \< \dt\Rho , \calN(\rho-\mz,\rhoG-\mz) >_{\calV} 
  \non
  \\
  && \quad {}
  + \iO \nabla\mu \cdot \nabla\calNO(\rho-\mz,\rhoG-\mz)
  + \iG \nablaG\muG \cdot \nablaG\calNG(\rho-\mz,\rhoG-\mz)
  \non
  \\
  && \quad {}
  + \iO |\nabla\rho|^2
  + \iG |\nablaG\rhoG|^2
  + \iO \beta(\rho) (\rho-\mz)
  + \iG \betaG(\rhoG) (\rhoG-\mz)
  \non
  \\
  && = \iO \rho u \cdot \nabla(\calNO(\rho-\rhoz,\rhoG-\rhoGz))
  + \iO \mu (\rho-\mz) 
  + \iG \muG (\rhoG-\mz)
  \non
  \\
  && \quad {}
  - \iO \pi(\rho) (\rho-\mz)
  - \iG \piG(\rhoG) (\rhoG-\mz).
  \non
\Esist
All of the integrals involving $\mu$ and $\muG$ cancel out by \eqref{defN}.
Now, we owe to~\eqref{trickMZ}, keep just the positive contribution on the \lhs,
and move the other terms to the \rhs.
Next, we account for the \Lip\ continuity of $\pi$ and~$\piG$,
the \Holder\ and Sobolev inequalities,
our assumption \eqref{hpu} on~$u$,
\eqref{primastima}, and~\eqref{secondastima}.
It then results that, \Aet,
\Bsist
  &&  \iO |\nabla\rho|^2
  + \iG |\nablaG\rhoG|^2
  + \delta_0 \iO |\beta(\rho)|
  + \delta_0 \iG |\betaG(\rhoG)|
  \non
  \\[2mm]
  && \leq \,\norma{\dt\Rho}_* \, \norma{\calN(\rho-\mz,\rhoG-\mz)}_{\calVz}
  + c \, \normaHH{(\pi(\rho),\piG(\rhoG))} \, \normaHH{(\rho-\mz,\rhoG-\mz)}
  \non
  \\[1mm]
  && \quad {}
  + \norma\rho_6 \, \norma u_3 \, \norma{\nabla(\calNO(\rho-\rhoz,\rhoG-\rhoGz))}_2
  \non
  \\[2mm]
  && \leq \,\norma{\dt\Rho}_* \, \norma{(\rho-\mz,\rhoG-\mz)}_*
  + c \, \bigl( \normaHH\Rho^2 + 1 \bigr)
  \non
  \\[1mm]
  && \quad {}
  + c \, \normaV\rho \, \norma{(\rho-\mz,\rhoG-\mz)}_*
  \,\leq\, c \,.
  \non
\Esist
We deduce (in particular) that
\Beq
  \zeta \in \LL\infty\Luno
  \aand
  \zetaG \in \LL\infty\LunoG .
  \non
\Eeq
Now, we test \eqref{seconda} by $(1,1)$
to obtain that, \Aet,
\Beq
  (|\Omega|+|\Gamma|) \mean\Mu 
  = \iO \zeta + \iG \zetaG
  + \iO \pi(\rho) + \iG \piG(\rhoG).
  \non
\Eeq
Thus, we can infer that
\Beq
  \mean\Mu \in L^\infty(0,+\infty) \,.
  \label{terzastima}
\Eeq

\medskip

It was already clear from \eqref{primastima} that the \omegalimit\ $\omega$ is nonempty.
Indeed, the weakly continuous $\calV$-valued function $\Rho$ is also bounded,
so that there exists a sequence $\tn\nearrow+\infty$ such that
the sequence $\graffe{\Rho(\tn)}$ is weakly convergent in~$\calV$.
More precisely, any sequence of times that tends to infinity contains a subsequence of this type.
Thus, it remains to prove the second part of the statement.
To this end, we fix an element $\Rhoo\in\omega$ 
and a corresponding sequence $\{\tn\}$ as in the definition~\eqref{omegalim}.
We also fix some \,$T\in(0,+\infty)$, set \aat
\Bsist
  && \mun(t) := \mu(\tn+t), \quad
  \rhon(t) := \rho(\tn+t), \quad
  \zetan(t) := \zeta(\tn+t),
  \non
  \\
  && \muGn(t) := \muG(\tn+t), \quad
  \rhoGn(t) := \rhoG(\tn+t), \quad
  \zetaGn(t) := \zetaG(\tn+t),
  \non
  \\
  && \un(t) := u(\tn+t),
  \non
\Esist
and notice that \eqref{hpu} implies that
\Beq
  \un \to 0
  \quad \hbox{strongly in $\L\infty{\Lx3}$}.
  \label{convun}
\Eeq
Moreover, it is clear that the 6-tuple $\soluzn$
satisfies the regularity conditions \Regsoluz\
and the equations \accorpa{prima}{terza} with $u$ replaced by~$\un$,
as well as the initial condition $\rhon(0)=\rho(\tn)$.
In particular, by construction, we have~that
\Beq
  (\rhon,\rhoGn)(0) \to \Rhoo
  \quad \hbox{weakly in $\calV$}.
  \label{limrhonz}
\Eeq
Furthermore, the global estimates already performed on $\soluz$
immediately imply some estimates on $\soluzn$ that are uniform with respect to~$n$.
Here is a list.
From \eqref{primastima} and \eqref{secondastima}, we infer that
\Beq
  \norma\Rhon_{\H1\calH\cap\L\infty\calV} \leq c \,.
  \label{stimarhon}
\Eeq
By virtue of \eqref{globint}, we also deduce that
\Bsist
  && (\nabla\mun,\nablaG\muGn) \to 0
  \quad \hbox{strongly in $(\L2\calH)^3$},
  \label{nablamun}
  \\[1mm]
  && (\dt\rhon,\dt\rhoGn) \to 0
  \quad \hbox{strongly in $\L2\calVzp$}.
  \label{dtrhon}
\Esist
On the other hand, \eqref{terzastima} yields a uniform estimate on the mean value \,$\mean\Mun$.
By combining this with \eqref{nablamun}, we conclude that
\Beq
  \norma\Mun_{\L2\calV} \leq \CT \,.
  \label{stimamun}
\Eeq
However, the estimates obtained till now are not sufficient to conclude,
and further estimates must be proved that ensure some better convergence for $\soluzn$ 
on the interval $(0,T)$.
Even though the argument is the same as in~\cite{CGS17},
we repeat it here for the reader's convenience,
at least in a short form.

\step 
First auxiliary estimate

We test \eqref{seconda} by the $\calV$-valued function $(\beta(\rhon),\beta(\rhoGn))$
and integrate with respect to time.
We have
\Bsist
  && \intQt \beta'(\rhon) |\nabla\rhon|^2
  + \intSt \beta'(\rhoGn) |\nablaG\rhoGn|^2
  + \intQt |\beta(\rhon)|^2
  + \intSt \betaG(\rhoGn) \, \beta(\rhoGn)
  \non
  \\
  && = \intQt \bigl( \mun - \pi(\rhon) \bigr) \, \beta(\rhon)
  + \intSt \bigl( \muGn - \piG(\rhoGn) \bigr) \, \beta(\rhoGn)\,,
  \label{perprimaaux}
\Esist
and we see that also the last integral on the \lhs\ is essentially \juerg{nonnegative. Indeed,} from \eqref{hpCC} and \eqref{hpBeta} it follows that
\Beq
  \intSt \betaG(\rhoGn) \, \beta(\rhoGn)
  \geq \frac 1 {2\eta} \intSt |\beta(\rhoGn)|^2  - \CT \,.
  \non
\Eeq
Thus, just the last integral on the \rhs\ needs some treatment.
We have
\Beq
  \intSt \muGn \, \beta(\rhoGn)
  \leq \frac 1 {4\eta} \intSt |\beta(\rhoGn)|^2
  + c \intSt |\muGn|^2
  \leq \frac 1 {4\eta} \intSt |\beta(\rhoGn)|^2
  + c \,,
  \non
\Eeq
thanks to \eqref{stimamun}.
By combining these inequalities, we conclude that
\Beq
  \norma\zetan_{\L2\HG} \leq \CT \,,
  \label{primaaux}
\Eeq
as well as an estimate for $\norma{\beta(\rhoGn)}_{\L2\HG}$, as a by-product.

\step
Second auxiliary estimate

We apply \cite[Lem.~3.1]{CGS13} \Aet\ to \eqref{seconda}, 
written for $\soluzn$, in the following form:
\Bsist
  && \iO \nabla\rhon \cdot \nabla v
  + \iG \nablaG\rhoGn \cdot \nablaG\vG
  + \iG \betaG(\rhoGn) \vG
  \non
  \\
  && = \iO \mun v 
  + \iG \muGn \vG  
  - \tauO \iO \dt\rhon \, v
  - \tauG \iG \dt\rhoGn \, \vG
  - \iO \bigl( \zetan + \pi(\rhon) \bigr) v
  - \iG \piG(\rhoGn) \vG .
  \non
\Esist
We obtain, in particular, that
\Beq
  \norma{\betaG(\rhoGn(t))}
  \leq c \, \bigl(
    \normaHH{\Mun(t)}
    + \normaH{\dt\Rhon(t)}
    + \normaH{(\zetan+\pi(\rhon))(t)}
  \bigr)
  \non
\Eeq
\Aat, where $c$ depends only on~$\Omega$.
By accounting for \eqref{stimarhon}, \eqref{stimamun}, and \eqref{primaaux}, 
we conclude that
\Beq
  \norma\zetaGn_{\L2\HG} \leq \CT \,.
  \label{secondaaux}
\Eeq

\step
Limits

Also this step closely follows~\cite{CGS17}.
Thanks to the estimates \eqref{stimarhon}, \eqref{stimamun}, \eqref{primaaux} and \eqref{secondaaux},
we have, for a subsequence which is still labeled by~$n$, 
\Bsist
  & \Rhon \to \Rhoi
  & \quad \hbox{weakly-star in $\H1\calH\cap\L\infty\calV$},
  \label{convRhon}
  \\
  & \Mun \to \Mui
  & \quad \hbox{weakly in $\L2\calV$},
  \label{convMun}
  \\
  & \Zetan \to \Zetai
  & \quad \hbox{weakly in $\L2\calH$}. 
  \label{convZetan}
\Esist
\juerg{The next step is to show} that $\soluzi$ satisfies \accorpa{prima}{terza} with \hbox{$u=0$.}
Namely, we can derive the integrated version
with time-dependent test functions $(v,\vG)\in\L2\calV$, as usual.
First of all, we notice that $\rhon\un$ converges to zero weakly in $\L\infty\Ldue$,
since $\rhon$ converges to $\rhoi$ weakly-star in $\L\infty{\Lx6}$ 
and \eqref{convun} holds.
Next, we derive from \eqref{convRhon} the strong convergence
\Beq
  \Rhon \to \Rhoi
  \quad \hbox{strongly in $\C0\calH$},
  \label{strongRhon}
\Eeq
by owing to the compact embedding $\calV\subset\calH$ 
and applying, e.g., \cite[Sect.~8, Cor.~4]{Simon}.
Hence, we can \juerg{identify} the limit of $(\pi(\rhon),\piG(\rhoGn))$ as $(\pi(\rhoi),\piG(\rhoGi))$
just by \Lip\ continuity.
This concludes the proof that \eqref{prima} and \eqref{seconda}
holds for the limiting 6-tuple in an integrated form,
which is equivalent to the pointwise formulation.
In order to derive that
$\zetai\in\beta(\rhoi)$ and $\zetaGi\in\betaG(\rhoGi)$,
we combine the weak convergence \eqref{convZetan} 
with the strong convergence \eqref{strongRhon}
and apply, e.g., \cite[Lemma~2.3, p.~38]{Barbu}.

\step
Conclusion

It remains to prove that the above limit leads to a stationary solution 
having the properties specified in the statement.
As in~\cite{CGS17}, we first derive that $\Rhoi$ belongs to $\L2\calW$ 
and solves the boundary value problem
\Bsist
  && - \Delta\rhoi + \zetai + \pi(\rhoi) = \mui
  \quad \aeQ\,,
  \label{bvpO}
  \\[1mm]
  && \dn\rhoi - \DeltaG\rhoGi + \zetaGi + \piG(\rhoGi) = \muGi
  \quad \aeS \,.
  \label{bvpG}
\Esist
Clearly, \eqref{dtrhon} implies that $\dt\Rhoi$ vanishes identically,
so that we are dealing with a time-dependent elliptic problem
in a variational form.
By using well-known estimates from trace theory
and from the theory of elliptic equations, and invoking a bootstrap argument
(see \cite{CGS17} for more details),
we successively draw the following conclusions:\\[1mm]
\hspace*{2.5mm}$\circ$ \,\,\eqref{bvpO} holds in the sense of distributions on~$\QT$;\\
\hspace*{2.5mm}$\circ$ \,\,$\Delta\rhoi\in\L2H$,
so that $\dn\rhoi$ is a well-defined element of $\L2{\HxG{-1/2}}$ 
satis-\linebreak \hspace*{7.7mm}fying the integration--by--parts formula in a generalized sense;\\
\hspace*{2.5mm}$\circ$ \,\,\eqref{bvpG} holds in a generalized sense;\\
\hspace*{2.5mm}$\circ$ \,\,$\DeltaG\rhoGi\in\L2{\HxG{-1/2}}$,
so that $\rhoGi\in\L2{\HxG{3/2}}$;\\
\hspace*{2.5mm}$\circ$ \,\,$\rhoi\in\L2W$;\\
\hspace*{2.5mm}$\circ$ \,\,$\dn\rhoi\in\L2\HG$, so that $\DeltaG\rhoGi\in\L2\HG$
and $\rhoGi\in\L2\WG$.

At this point, we can conclude the proof
by repeating the argument of \cite{CGS17} for the reader's convenience.
Since \gianni{both $\dt\Rhoi$ and $(\nabla\mui,\nablaG\muGi)$ vanish by~\accorpa{nablamun}{dtrhon}}, 
there exist $\Rhos\in\calV$ and $\musi\in L^2(0,T)$ such that
\Beq
  \Rhoi(x,t) = \Rhos(x)
  \aand
  \Mui(x,t) = (\musi(t),\musi(t))
  \quad \hbox{for a.a.\ $(x,t)\in\QT$} \,.
  \non
\Eeq
We show that $\Zetai$ is time independent as well and that $\musi$ is a constant
by accounting for~\eqref{perunicita}.
Assume that $\beta$ is single-valued.
Then, $\zetai=\beta(\rhoi)$ takes the value $\zetas:=\beta(\rhos)$ at any time.
Therefore, \eqref{bvpO} implies that $\mui$ is time independent as well,
so that the function $\musi$ is a constant that we term~$\mus$.
Thus, the \rhs\ of \eqref{bvpG} is the same constant~$\mus$.
As this does not depend on time, the same holds for~$\zetaG$,
which takes some value $\zetaGs\in\HG$ \aet.
Assume now that $\betaG$ is single-valued.
Then, we first use \eqref{bvpG} to derive that 
$\zetaGi=\betaG(\rhoGi)$ and $\muGi$ are time independent.
In particular, $\musi$~attains some constant value~$\mus$,
so that $\zetai$ is time independent, by comparison in~\eqref{bvpO}.
\juerg{Thus,  in both cases}
the quadruplet $\soluzs$ is a stationary solution
corresponding to the value $\mus$ of the chemical potential.
Finally, we have that $\Rhos=\Rhoo$.
Indeed, \eqref{convRhon} implies weak convergence also in $\C0\calH$, whence
\Beq
  \Rhon(0) \to \Rhoi(0) = \Rhos
  \quad \hbox{weakly in~$\calH$},
  \non
\Eeq
and we can compare with \eqref{limrhonz}.\QED


\section{The optimal control problem}
\label{CONTROL}
\setcounter{equation}{0}

In this section, we are going to prove the Theorems~\ref{Optimum} and~\ref{Necessary}.
Thus, we fix $T>0$ and consider all problems on the finite time interval~$(0,T)$.
We adopt the ideas of~\cite{CGS16}
and take the problem \Pbl, written with $\tauO=\tauG=\tau\in(0,1)$, 
as an approximating problem for the pure case. We thus 
work with problems of the form
\Bsist
  && \< \dt\Rhot , (v,\vG) >_{\calV}
  - \iO \rhot u \cdot \nabla v
  + \iO \nabla\mut \cdot \nabla v
  + \iG \nablaG\muGt \cdot \nablaG\vG
  = 0
  \non
  \\[1mm]
  && \quad \hbox{\aet\ and for every $(v,\vG)\in\calV$},
  \label{primatau}
  \\
  \separa
  && \tau \iO \dt\rhot \, v
  + \tau \iG \dt\rhoGt \, \vG
  + \iO \nabla\rhot \cdot \nabla v
  + \iG \nablaG\rhoGt \cdot \nablaG\vG
  \non
  \\[2mm]
  && \quad {}
  + \iO (\beta+\pi)(\rhot) \bigr) v
  + \iG (\betaG+\piG)(\rhoGt) \vG
  = \iO \mut v 
  + \iG \muGt \vG
  \non
  \\[1mm]
  && \quad \hbox{\aet\ and for every $(v,\vG)\in\calV$},
  \label{secondatau}
  \\[2mm]
  && \rhot(0) = \rhoz, 
  \label{cauchytau}
\Esist
for given $u\in\Uad$ and \juerg{$\tau\in (0,T)$. Let us assume that
the assumptions  \eqref{hptau}--\eqref{perunicita}, \eqref{hprhoz}--\eqref{hpmz}, and \eqref{hprhozpure}
are fulfilled. Then we can infer from \cite[Thm.~2.6]{CGS13} that
the problem \eqref{primatau}--\eqref{cauchytau} has for every $u\in\Uad$ and every $\tau\in (0,T)$
a unique solution $(\mut,\muGt,\rhot,\rhoGt)$ such that}
\begin{align}
\label{regt}
&\juerg{(\mut,\muGt)\in L^\infty(0,T;\calW),\quad 
(\rhot,\rhoGt)\in W^{1,\infty}(0,T;\calH)\cap H^1(0,T;\calV)\cap L^\infty(0,T;\calW).}
\end{align}
\juerg{In particular, we have that $\mut,\rhot\in L^\infty(Q)$ and $\muGt,\rhoGt\in L^\infty(\Sigma)$.}

\juerg{We now aim to derive bounds for the solutions that are uniform with respect to $\tau\in (0,1)$ and 
$u\in\Uad$. In particular, we establish a uniform $L^\infty$ estimate for the solution to \Pbl.
This was already announced in the Remarks~2.7 and~7.1 of~\cite{CGS13},
but not proved, essentially.
We sketch the derivation here, for the reader's convenience.
For some of the steps, we proceed very quickly
and refer for the details to the proofs of the estimates of~\cite[Sect.~6]{CGS13}.
It the following, the symbol $c$ denotes constants that are independent of $\tau\in(0,1)$ and $u\in\Uad$.
These constants may however depend on the fixed final time~$T$,
although we simply write~$c$.}

\step
Uniform estimates

Suppose that $u\in\Uad$.
The argument used to prove the estimate \eqref{primastima} for the solution to \Pbl\
applies to the system \accorpa{primatau}{cauchytau} as well and yields
(see also \eqref{globint} and~\eqref{stimef})
\begin{align}
  & \norma{\nabla\mut}_{(\L2H)^3}
  + \norma{\nablaG\muGt}_{(\L2\HG)^3}
  + \norma\Rhot_{\H1\calVp\cap\L\infty\calV}
  \leq c\,,
  \qquad \label{unit1}
    \\[1mm]
   & \norma{f(\rhot)}_{\L\infty\Luno}
  + \norma{\fG(\rhoGt)}_{\L\infty\LunoG}
  \leq c \,.\label{unit2}
  \end{align}
Then, we can follow the procedure used to derive \eqref{secondastima}
and obtain
\Beq
\label{unit3}
  \norma{\dt\Rhot}_{\L\infty\calVp} \leq c\,, 
 \Eeq
since the strong assumption \eqref{hprhozpure} on~$\rhoz$ is assumed to hold true.
Now, we define \,$\hm(t):=\mean\Mut(t)$\, \aat\
and test \eqref{primatau}, written at the time~$t$, by
$(\mut(t)-\hm(t),\muGt(t)-\hm(t))\in \calVz$.
We obtain, \aet,
\Beq
  \iO |\nabla\mut|^2
  + \iG |\nablaG\muGt|^2
  = - \< \dt\Rhot , (\mut-\hm,\muGt-\hm) >_{\calV}
  + \iO \rhot u \cdot \nabla\mut .
  \non
\Eeq
By using the norms \eqref{normaVz} and~\eqref{normastar},
as well as Young's inequality and the above estimates,
we deduce that
\Bsist
  && \iO |\nabla\mut|^2
  + \iG |\nablaG\muGt|^2
  \leq \norma{\dt\Rhot}_* \, \norma{(\mut-\hm,\muGt-\hm)}_{\calVz}
  + \juerg{\norma u_\infty \, \norma\rhot_2 \, \|\nabla\mut\|_2}
  \non
  \\
  && \leq \frac 34 \iO |\nabla\mut|^2
  + \frac 12 \iG |\nablaG\muGt|^2
  + c \,,
  \non
\Esist
\gianni{so that we can infer that
\Beq
  \norma{(\nabla\mut,\nablaG\muGt)}_{(\L\infty\calH)^3}
  \leq c, 
  \quad \hbox{whence} \quad
  \norma{(\mut-\hm,\muGt-\hm)}_{\L\infty\calH}
  \leq c \,.
  \non
\Eeq
}%
Then, we can come back to the procedure adopted to \gianni{show} the validity of \eqref{terzastima}.
It is clear that we can deduce that
the mean value $\hm$ is bounded in~$L^\infty(0,T)$.
We thus conclude that
\begin{equation}
\label{unit4}
  \norma\Mut_{\L\infty\calV} \leq c \,.
  \end{equation}
In addition, since \eqref{hpCC} implies that
\Beq
  \betaG(r) \beta(r)
  \geq \frac 1 {2\eta} \, |\beta(r)|^2 - c
  \quad \hbox{for every $r\in\erre$}
  \non
\Eeq
(in fact, \juerg{this also holds true for the Yosida regularizations 
$\betaeps$ and~$\betaGeps$, 
as it was} proved in~\cite[Lemma 4.4]{CaCo}),
we can test \eqref{seconda} by $(\beta(\rhot),\beta(\rhoGt))$
(without integrating in time) and deduce~that
\Beq
  \norma{\beta(\rhot)}_{\L2H} + \norma{\beta(\rhoGt)}_{\L2\HG} \leq c \,.
  \non
\Eeq
Therefore, we can use the parts $vi)$ and $v)$ of \cite[Thm.~2.1]{CGS13}
(indeed, this procedure should be performed on the $\eps$-approximating problem
whose nonlinearities are \Lip\ continuous
and \juerg{thus} yield constants that do not depend on~$\eps$)
and infer that
\Beq
  \norma{\betaG(\rhoGt)}_{\L\infty\HG} \leq c
  \aand
  \norma\Rhot_{\L\infty\calW} \leq c \,.
  \label{unit5}
\Eeq
From the latter we conclude that
\Beq
  \norma\Rhot_\infty \leq \hR\,,
  \label{stimaunif}
\Eeq
where the constant \,$\hR$\, which we have marked with a special symbol 
depends only on the structure of the system and of the control box $\Uad$,
as well as on the norms involved in our assumptions on the initial datum.
Hence, this estimate is uniform with respect to \juerg{$\tau\in (0,T)$ and
$u\in\Uad$}.

\juerg{We draw some important consequences from \eqref{stimaunif}. To this end,
we introduce for \,$\tau\geq0$\, the solution} operators 
\Bsist
  && \calS_\tau : \Uad \to L^2(0,T;\calV) \times \bigl( H^1(0,T;\calVp)
   \cap \L\infty\calV \bigr),
  \non
  \\
  && \calS^2_\tau : \Uad \to H^1(0,T;\calVp)
   \cap L^\infty(0,T;\calV),
     \label{defcalS}
\Esist
which are defined as follows:
for a given $u\in\Uad$, the value $\calS_\tau(u)$ is the solution $(\mut,\muGt,\rhot,\rhoGt)$
to the system \accorpa{primatau}{cauchytau} 
given by Theorem~\ref{RecallCGS},
and $\calS^2_\tau(u)$
is its second component pair~$\Rhot$. 

\juerg{Now let $\tau\in (0,1)$ and $u\in\Uad$ be arbitrary. Then, with $(\rhot,\rhoGt)={\cal S}_\tau^2(u)$,
it follows from \eqref{stimaunif} that}
\Beq
  \norma{f^{(j)}(\rhot)}_\infty
  + \norma{\fG^{(j)}(\rhoGt)}_\infty
  \leq c
  \quad \hbox{for $0\leq j\leq3$},
  \label{dastimaunif}
\Eeq
since the potentials are smooth. 
In particular, the functions 
\Beq
  \psit := f''(\rhot)
  \aand 
  \psiGt := \fG''(\rhoGt) ,
  \label{defpsit}
\Eeq
are \juerg{bounded in $L^\infty(Q)$ and $L^\infty(\Gamma)$, respectively, uniformly with respect to $\tau\in (0,T)$ and $u\in\Uad$. Moreover, the functions}
\gianni{%
\Bsist
  && \phi_3^\tau:=\beta_3(\rhot-\widehat\rho_Q), \quad
  \phi_4^\tau:=\beta_4(\rhoGt-\widehat\rho_\Sigma),
  \non
  \\
  && \phi_5^\tau:=\beta_5(\rhot(T)-\widehat\rho_\Omega),
  \aand
  \phi_6^\tau:=\beta_6(\rhoGt(T)-\widehat\rho_\Gamma),
  \label{defphit}
\Esist
}%
\juerg{are bounded in $L^2(Q)$, $L^2(\Sigma)$, $L^2(\Omega)$, and $L^2(\Gamma)$, respectively, uniformly
with respect to $\tau\in (0,T)$ and $u\in\Uad$.}
We also remark that
\Beq
  \hbox{the functions $f^{(j)}$ and $\fG^{(j)}$ are \Lip\ continuous on $[-\hR,\hR]$ for $0\leq j\leq2$}.
  \label{lipf}
\Eeq

\vspace{3mm}
\juerg{After these preparations, we are now ready to go on with our project for the proof of our results.
The following approximation result resembles Theorem~\cite[Thm.~3.1]{CGS16}, which was 
established for the case of logarithmic potentials and $\tauO>0$, $ \tauG>0$.}

\Bthm
\label{Convergence}
\juerg{
Suppose that the assumptions  \eqref{hpBeta}--\eqref{hpCC}, \eqref{hprhoz}, \eqref{hprhozpure}, and 
\eqref{hppotcon}--\eqref{hpUad},
are fulfilled, and assume that sequences $\{\tau_n\}\subset (0,1)$ and $\{\utn\}\subset\Uad$ are given such that
$\tau_n\searrow 0$ and $\utn\to u$ weakly-star in $\calX$ for some $u\in\Uad$. Then, with
$(\mutn,\muGtn,\rhotn,\rhoGtn)=\calS_{\tau_n}(\utn)$ and $(\mu,\muG,\rho,\rhoG)=\calS_0(u)$, we have that
}
\begin{align}
\label{conmu}
  &\juerg{ \Mutn \to \Mu
  \quad \hbox{weakly-star in \,$L^\infty(0,T;\calV)$}},
  \\[1mm]
  \label{conrho}
  &\juerg{ \Rhotn \to \Rho
  \quad \hbox{weakly-star in \,$W^{1,\infty}(0,T;\calV^*)\cap L^\infty(0,T;\calW)$}}\non\\
  &\juerg{\hspace*{38.2mm}\mbox{and strongly in \,$C^0(\overline{Q})\times C^0(\overline\Sigma) $}.}
 \end{align}
\juerg{Moreover, it holds that}
\begin{align}
& \calJ(\calS_0^2(u),u)
 \, \leq\, \liminf_{n\to\infty} \calJ(\Rhotn,\utn)
  \label{conj1}
  \\
& \calJ(\calS^2_0(v),v)
 \, \,=\, \lim_{n\to\infty} \calJ(\calS^2_{\tau_n}(v),v)
  \quad \hbox{for every $v\in\Uad$}.
  \label{conj2}
\end{align}
\Ethm

\vspace{2mm}\noindent
{\em Proof:} \quad \juerg{
\,Let $\,\{\tau_n\}\subset (0,1)\,$ be any sequence such that $\tau_n\searrow 0$ as $n\to\infty$, and suppose that 
$\{\utn\}\subset\Uad$ converges weakly-star in ${\cal X}$ to some $u\in\Uad$. By virtue of the global estimates \eqref{unit3}--\eqref{unit5}, there are some subsequence of $\{\tau_n\}$, which is again indexed by $n$, and two pairs $\,\Mu,\Rho\,$ such that
\eqref{conmu} and the first convergence result of \eqref{conrho} hold true. It then follows from 
standard compact embedding results (cf. \cite[Sect.~8, Cor.~4]{Simon}) that} 
\begin{align}
\label{strcon1}
&\juerg{\rho^{\tau_n}\to\rho\quad\mbox{strongly in }\,L^2(0,T;V)\cap C^0(\overline Q),}
\end{align}
\juerg{which also implies that}
\begin{equation}
\label{strcon2}
\juerg{\rho_\Gamma^{\tau_n}\to\rhoG\quad\mbox{strongly in }\,C^0(\overline \Sigma)\,.}
\end{equation}
\juerg{In particular, $(\rho(0),\rhoG(0))=(\rho_0,\rho_{0|\Gamma})$ and $\,\rhoG=\rho_{|\Sigma}$. In addition, 
we obviously have that}
\begin{align}
\label{strcon3}
& \juerg{\beta(\rho^{\tau_n})\to \beta(\rho), \quad {\pi}(\rho^{\tau_n})\to {\pi}(\rho),\quad\mbox{both strongly in $\,C^0(\overline Q)$},}\\[1mm]
\label{strcon4}
&\juerg{\betaG(\rho^{\tau_n}_\Gamma)\to\betaG(\rhoG), \quad {\piG}(\rho_\Gamma^{\tau_n})\to{\piG}(\rhoG),\quad\mbox{both
strongly in $\,C^0(\overline \Sigma)$}.}
\end{align} 
\juerg{Moreover, it is easily verified that, at least weakly in $L^1(Q)$,}
\begin{align}
&\juerg{\nabla\rho^{\tau_n}\cdot\utn\to \nabla\rho\cdot u\,.}
\end{align}
\juerg{Combining the above convergence results, we may pass to the limit
as $n\to\infty$ in the equations \eqref{primatau}--\eqref{cauchytau} (written for 
$\tau=\tau_n$ and $u=\utn$) to find that $(\mu,\muG,\rho,\rhoG)$ and $\,u\,$
satisfy the equations \eqref{prima}--\eqref{cauchy} for $\tauO=\tauG=0$ with
$\zeta=\beta(\rho)$ and $\zeta_\Gamma=\betaG(\rhoG)$. Owing to the uniqueness result
stated in Theorem 2.1, we therefore have that
 $(\mu,\muG,\rho,\rhoG)=\calS_0(u)$, and since the limit is unique, the convergence 
properties \eqref{conmu} and \eqref{conrho} hold true for the entire
sequences.}

\juerg{It remains to show the validity of \eqref{conj1} and \eqref{conj2}. 
In view of \eqref{conrho}, the inequality \eqref{conj1}  is an immediate consequence of the 
weak and weak-star sequential semicontinuity properties of the cost functional ${\cal J}$. To establish
the identity \eqref{conj2}, let $v\in \Uad$ be arbitrary  and  put
$(\rho^{\tau_n},\rhoG^{\tau_n})={\cal S}_{\tau_n}^2(v)$, for $n\in\enne$. 
Taking Theorem 2.1 into account, and arguing as in the first part of this proof,
we can conclude that $\{{\calS}_{\tau_n}^2(v)\}$ converges to $\Rho={\calS}_0^2(v)$ in the sense
of \eqref{conrho}. In particular, }
$$
\juerg{{\calS}_{\tau_n}^2(v)\to {\cal S}_0^2(v)\quad\mbox{strongly in }\,C^0(\overline Q)\times C^0(\overline\Sigma).}
$$
\juerg{As the cost functional ${\cal J}$ is obviously continuous in the variables $(\rho,\rhoG)$
with respect to the strong topology of $\,C^0(\overline Q)\times C^0(\overline\Sigma)$,
we thus conclude that \eqref{conj2} is valid.} \QED

\vspace{2mm}
As a  corollary of Theorem~\ref{Convergence}, we can prove Theorem~\ref{Optimum}.

\vspace{2mm}\noindent
\juerg{{\em Proof of Theorem~\ref{Optimum}}: \quad At first, we observe that all the assumptions
for an application of the arguments employed in the proof of Theorem 4.1 are fulfilled. We now conclude
in two steps.}

\noindent
\underline{{\sc Step 1:} }

\vspace{1mm}\noindent
\juerg{We first consider the problem of minimizing the cost functional $\,{\cal J}\,$ subject to $u\in\Uad$
and to the state system \eqref{prima}--\eqref{cauchy} for fixed $\tauO=\tauG=\tau>0$. We claim that this
optimal control problem, which we denote by $({\cal P}_\tau)$, admits at least one optimal pair
for every $\tau>0$. Indeed, let $\tau>0$ be fixed, and let $((\mu_n,\mu_{n_\Gamma},\rho_n,\rho_{n_\Gamma}),
u_n)$, $n\in\enne$, be a minimizing sequence for $({\cal P}_\tau)$, that is, we assume that we have 
\,$(\mu_n,\mu_{n_\Gamma},\rho_n,\rho_{n_\Gamma})={\cal S}_\tau(u_n)$ for all $n\in\enne$ and} 
$$
\juerg{
\lim_{n\to\infty}\,{\cal J}((\rho_n,\rho_{n_\Gamma}),u_n)\,=\,\inf_{v\in\Uad}\,J({\cal S}^2_\tau(v),v)
\,=:\,\sigma\,\ge\,0\,.}
$$
\juerg{Then, by the same token as in the proof of Theorem 4.1, there are $(\mu,\muG,\rho,\rhoG)$ and 
$u\in\Uad$ with $\,(\mu,\muG,\rho,\rhoG)={\cal S}_\tau(u)\,$ such that (cf. \eqref{conmu} and \eqref{conrho}) }
\begin{align*}
&\juerg{(\mu_n,\mu_{n_\Gamma})\to (\mu,\muG)\quad\mbox{weakly-star in $L^\infty(0,T;\calV)$,}  }\\
&\juerg{(\rho_n,\rho_{n_\Gamma})\to (\rho,\rhoG) \quad\mbox{strongly in $C^0(\overline Q)\times
C^0(\overline\Sigma)$.}}
\end{align*}
\juerg{The sequential lower semicontinuity properties of ${\cal J}$ then imply that 
$((\mu,\mu_\Gamma,\rho,\rhoG),u)$ is optimal for $({\cal P}_\tau)$, which proves the claim.}

\vspace{2mm}\noindent
\juerg{\underline{\sc Step 2:}} 

\vspace{1mm}\noindent
\juerg{We now pick an arbitrary sequence $\{\tau_n\}$ such that $\tau_n\searrow0$ as $n\to\infty$.  
Then, as has been shown in Step 1, the optimal control problem (${\cal P}_{\tau_n}$) has for
every $n\in\enne$ an optimal pair $((\mu^{\tau_n},\mu_\Gamma^{\tau_n},\rho^{\tau_n},\rho_\Gamma^{\tau_n}),\utn)$, where $\,\utn\in\Uad\,$ and  $\,(\rho^{\tau_n},\rho_\Gamma^{\tau_n})={\cal S}^2_{\tau_n}(\utn)$. 
Since $\Uad$ is a bounded subset of ${\cal X}$, we may without loss of generality assume
that $\,\utn\to u\,$ weakly-star in ${\cal X}$ for some $\,u\in\Uad$. From Theorem 4.1 we infer that
with $(\mu,\muG,\rho,\rhoG)={\cal S}_0(u)$ the convergence properties \eqref{conmu} and \eqref{conrho}
are valid, as well as \eqref{conj1}.
Invoking the optimality of $((\mu^{\tau_n},\mu_\Gamma^{\tau_n},\rho^{\tau_n},\rho_\Gamma^{\tau_n}),\utn)$ for (${\cal P}_{\tau_n}$) and \eqref{conj2},
we then find, for every $v\in\Uad$, that} 
\begin{align}
\label{tr3.3}
&\juerg{{\cal J}(\Rho,u)\,=\,{\cal J}({\cal S}_{0}^2(u),u)\,\le\,
\liminf_{n\to\infty}\,{\cal J}({\mathcal S}_{\tau_n}^2(u^{\tau_n}),u^{\tau_n})}
 \nonumber\\[1mm]
&\juerg{\leq\,\liminf_{n\to\infty}\,{\cal J}({\cal S}_{\tau_n}^2(v),v)\, =\,\lim_{n\to\infty} {\cal J}({\cal S}_{\tau_n}^2(v),v)\,=\,
{\cal J}({\mathcal S}_{0}^2(v),v),}
\end{align}  
\juerg{which yields that $u$ is an optimal control for the control problem for $\tau=0$ with the associate state
$(\mu,\muG,\rho,\rhoG)$. The assertion is thus proved.
\QED}

\vspace{5mm}
\juerg{We remark at this place that the existence result of Theorem 2.6 can also be proved
using the same direct argument as in Step 1 of the above proof. We have chosen to employ Theorem 4.1, here, since
it shows that, for small $\tau>0$, optimal pairs for $({\cal P}_\tau)$ are likely to be
`close' to optimal pairs for the case $\tau=0$. However, the result does not yield any information on whether every solution to the optimal control problem for $\,\tau=0\,$ can be approximated by a sequence of solutions to the problems $({\mathcal{P}}_{\tau_n})$. Unfortunately, we are not able to prove such a 
general `global' result. Instead, we 
can only give a `local' answer for every individual optimizer of the control problem for $\tau=0$. For this purpose,
we employ a trick due to Barbu~\cite{Barbu}. To this end, let $\bar u\in\Uad$
be an arbitrary optimal control for $\tau=0$, and let $(\bar\mu,\bar\mu_\Gamma,
\bar\rho,\bar\rho_\Gamma)$ be the associated solution  to the state system (\ref{prima})--(\ref{cauchy}) 
for $\tauO=\tauG=0$. In particular, $\,(\bar\rho,\bar\rho_\Gamma)={\cal S}_0^2 (\bar u)$. We associate with this 
optimal control the {\em adapted cost functional}}
\begin{equation}
\label{cost2}
\juerg{\widetilde{\cal J}((\rho,\rhoG),u):={\cal J}((\rho,\rhoG),u)\,+\,\frac{1}{2}\,\|u-\bar{u}\|^2_{(L^2(Q))^3}}
\end{equation}
\juerg{and, for every $\tau\in (0,1)$, a corresponding {\em adapted optimal control problem},}

\vspace{4mm}\noindent
\juerg{($\widetilde{\mathcal{P}}_{\tau}$)\quad Minimize $\,\, \widetilde {\cal J}(\Rho,u)\,\,$
for $\,u\in\Uad$, subject to the condition that  
(\ref{primatau})--(\ref{cauchytau}) \hspace*{11mm} be satisfied.}

\vspace{3mm}
\juerg{With the same direct argument as in Step 1 of the proof of Theorem 2.6 (which needs no repetition, here), we 
can show that under the assumptions of Theorem 2.6 the optimal control problem 
$(\widetilde{\cal P}_\tau)$ admits for every $\tau\in (0,1)$ a solution.
The next result is an analogue of Theorem~\cite[Thm.~3.4]{CGS16}, which has been shown for potentials
of logarithmic type and $\tauO>0$ and $\tauG>0$:}

\Bthm
\label{Strongconv}
\juerg{Suppose that the assumptions  \eqref{hpBeta}--\eqref{hpCC}, \eqref{hprhoz}, \eqref{hprhozpure}, and 
\eqref{hppotcon}--\eqref{hpUad},
are fulfilled. Moreover, let $\ub\in\Uad$ be an optimal control related to the cost functional \eqref{cost}
and to the state system \Pbl\ for $\tauO=\tauG=0$, where $(\mub,\muGb,\rhob,\rhoGb)={\cal S}_0(\ub)$ is the corresponding  state. If  $\{\tau_n\}\subset(0,1)$ is any sequence with $\tau_n\searrow0$ as $n\to\infty$,
then there exist a subsequence $\{\tau_{n_k}\}$
and, for every~$k$, an optimal control $\utnk\in\Uad$ of the adapted control problem
$(\widetilde{\cal P}_{\tau_{n_k}})$, such that, as~$k\to\infty$,}
\Beq
 \juerg{ \utnk \to \ub
  \quad \hbox{strongly in $(\LQ2)^3$},}
  \label{strongconv}
\Eeq
\juerg{and such that the convergence properties \accorpa{conmu}{conrho} are satisfied,
where $\{\tau_n\}$ and $(\mu,\muG,\linebreak \rho,\rhoG)$ are replaced by $\{\tau_{n_k}\}$ and $(\mub,\muGb,\rhob,\rhoGb)$,
respectively.
Moreover, we have}
\Beq
  \lim_{k\to\infty} \tilde\calJ((\rhotnk,\rhoGtnk),\utnk)
  = \calJ((\rhob,\rhoGb),\ub).
  \label{convJ}
\Eeq
\Ethm

\vspace{2mm}\noindent
{\em Proof:} \quad \juerg{
Let $\tau_n \searrow 0$ as $n\to\infty$. For any $ n\in\enne$, we pick an optimal control 
$u^{\tau_n} \in \Uad\,$ for the adapted problem $(\widetilde{\cal P}_{\tau_n})$ and denote by 
$(\mu^{\tau_n},\mu_\Gamma^{\tau_n},\rho^{\tau_n},\rho_\Gamma^{\tau_n})$ the associated solution to the problem (\ref{primatau})--(\ref{cauchytau}) for $\tau=\tau_n$ and $u=\utn$. 
By the boundedness of $\Uad$ in $\calX$, there is some subsequence $\{\tau_{n_k}\}$ of $\{\tau_n\}$ such that}
\begin{equation}
\label{ugam}
\juerg{u^{\tau_{n_k}}\to u\quad\mbox{weakly-star in }\,{\cal X}
\quad\mbox{as }\,k\to\infty,}
\end{equation}
\juerg{with some $u\in\Uad$. Thanks to Theorem 4.1, the convergence properties \eqref{conmu}--\eqref{conrho}
hold true, where $(\mu,\muG,\rho,\rhoG)$ is the unique solution to the state system
\eqref{prima}--\eqref{cauchy} for $\tauO=\tauG=0$. In particular,  $((\rho,\rhoG),u)$
is admissible for the non-adapted control problem with the cost functional \eqref{cost}.}

\vspace{2mm}\quad
\juerg{We now aim to prove that $u=\bar u$. Once this is shown, then the uniqueness result of 
Theorem 2.1 yields that also $(\mu,\muG,\rho,\rhoG)=(\mub,\muGb,\rhob,\rhoGb)$, 
which implies that the properties \eqref{conmu}--\eqref{conrho} are satisfied, where $(\mu,\muG,\rho,
\rhoG)$ is replaced by  $(\bar\mu,\bar\mu_\Gamma,\bar\rho,\bar\rho_\Gamma)$.}

\juerg{Now observe that, owing to the weak sequential lower semicontinuity of 
$\widetilde {\cal J}$, 
and in view of the optimality property of $((\bar\rho,\bar\rho_\Gamma),\ub)$,}
\begin{align}
\label{tr3.6}
&\juerg{\liminf_{k\to\infty}\, \widetilde{\cal J}((\rho^{\tau_{n_k}},
\rho_\Gamma^{\tau_{n_k}}),u^{\tau_{n_k}})
\ge \,{\cal J}((\rho,\rhoG),u)\,+\,\frac{1}{2}\,
\|u-\bar{u}\|^2_{(L^2(Q))^3}}\nonumber\\[1mm]
&\juerg{\geq \, {\cal J}((\bar\rho,\bar\rho_\Gamma),\bar u)\,+\,\frac{1}{2}\,\|u-\bar{u}\|^2_{(L^2(Q))^3}\,.}
\end{align}
\juerg{On the other hand, the optimality property of  $\,((\rho^{\tau_{n_k}},\rho_\Gamma^{\tau_{n_k}}),u^{\tau_{n_k}})
\,$ for problem $(\widetilde {\cal P}_{\tau_{n_k}})$ yields that
for any $k\in\enne$ we have}
\begin{equation}
\label{tr3.7}
\juerg{\widetilde {\cal J}((\rho^{\tau_{n_k}},\rho_\Gamma^{\tau_{n_k}}),u^{\tau_{n_k}})\, =\,
\widetilde {\cal J}({\cal S}^2_{\tau_{n_k}}(u^{\tau_{n_k}}),
u^{\tau_{n_k}})\,\le\,\widetilde {\cal J}({\cal S}^2_{\tau_{n_k}}
(\bar u),\bar u)\,,}
\end{equation}
\juerg{whence, taking the limit superior as $k\to\infty$ on both sides and invoking (\ref{conj2}) in
Theorem~4.1,}
\begin{align}
\label{tr3.8}
&\juerg{\limsup_{k\to\infty}\,\widetilde {\cal J}((\rho^{\tau_{n_k}},
\rho_\Gamma^{\tau_{n_k}}),u^{\tau_{n_k}}) }
\nonumber\\[1mm]
&\juerg{\le\,\widetilde {\cal J}({\calS}_0^2(\bar u),\bar u) 
\,=\,\widetilde {\cal J}((\bar\rho,\bar\rho_\Gamma),\bar u)
\,=\,{\cal J}((\bar\rho,\bar\rho_\Gamma),\bar u)\,.}
\end{align}
\juerg{Combining (\ref{tr3.6}) with (\ref{tr3.8}), we have thus shown that 
$\,\frac{1}{2}\,\|u-\bar{u}\|^2_{(L^2(Q))^3}=0$\,,
so that $\,u=\bar u\,$  and thus also $\,(\mu,\muG,\rho,\rhoG)
=(\bar\mu,\bar\mu_\Gamma,\bar\rho,\bar\rho_\Gamma)$. 
Moreover, (\ref{tr3.6}) and (\ref{tr3.8}) also imply that}
\begin{align}
\label{tr3.9}
&\juerg{
{\cal J}((\bar\rho,\bar\rho_\Gamma),\bar u) \, =\,\widetilde{\cal J}((\bar\rho,\bar\rho_\Gamma),\bar u)
\,=\,\liminf_{k\to\infty}\, \widetilde{\cal J}((\rho^{\tau_{n_k}},
\rho_\Gamma^{\tau_{n_k}}), u^{\tau_{n_k}}) }\nonumber\\[1mm]
&\juerg{\,=\,\limsup_{k\to\infty}\,\widetilde{\cal J}((\rho^{\tau_{n_k}},
\rho_\Gamma^{\tau_{n_k}}), u^{\tau_{n_k}}) \,
=\,\lim_{k\to\infty}\, \widetilde{\cal J}((\rho^{\tau_{n_k}},
\rho_\Gamma^{\tau_{n_k}}), u^{\tau_{n_k}})\,,}
\end{align}                                     
\juerg{which proves {(\ref{strongconv})} and, at the same time, also (\ref{convJ}). This concludes the proof
of the assertion.\QED}

\vspace{5mm}
In order to prove Theorem~\ref{Necessary},
we need estimates that are uniform with respect to $\tau\in (0,1)$
for the solution $(\qt,\qGt,\pt,\pGt)$ 
to the approximating adjoint system
corresponding to a velocity field $\juerg{\bar u^\tau}\in\Uad$
and to the associated solution $(\mut,\muGt,\rhot,\rhoGt)\juerg{={\cal S}_\tau(\bar u^\tau)}$
to the state system \accorpa{primatau}{cauchytau} with $u=\bar u^\tau$.
\juerg{The approximating adjoint} system reads
\Bsist
  && - \< \dt (\pt+\tau\qt,\pGt+\tau\qGt) , (v,\vG) >_{\calV}
  + \iO \nabla\qt \cdot \nabla v
  + \iG \nablaG\qGt \cdot \nablaG\vG
  \non
  \\
  && \quad {}
  + \iO \psit\qt v
  + \iG \psiGt\qGt \vG
  - \iO \ubt \cdot \nabla \juerg{\pt} \, v
  = \iO \juerg{\phi_3^\tau} \, v
  + \iG \juerg{\phi_4^\tau} \, \vG
  \non
  \\[1mm]
  && \quad \hbox{\aet\ and for every $(v,\vG)\in\calV$},
  \label{primaAtau}
  \\
  \separa
  && \iO \nabla\pt \cdot \nabla v
  + \iG \nablaG\pGt \cdot \nablaG\vG
  = \iO \qt v 
  + \iG \qGt \vG
  \non
  \\[1mm]
  && \quad \hbox{\aet\ and for every $(v,\vG)\in\calV$},
  \label{secondaAtau}
  \\
  && \< (\pt+\tau\qt,\pGt+\tau\qGt)(T) , (v,\vG) >_{\calV}
  = \iO \juerg{\phi_5^\tau} \, v
  + \iG \juerg{\phi_6^\tau} \,\vG
  \non
  \\
  && \quad \hbox{for every $(v,\vG)\in\calV$}, 
  \label{cauchyAtau}
\Esist
\juerg{with the functions defined in \eqref{defpsit} and \eqref{defphit}}.
The velocity field $\ubt$ appearing in the above system
could be any element of~$\Uad$, in principle.
However, we will use \accorpa{primaAtau}{cauchyAtau} only when this field is
an optimal control for the control problem associated with the adapted cost functional \eqref{cost2}
and to the $\tau$-state system \juerg{\accorpa{primatau}{cauchytau}}.

The basic idea is the following:
we integrate equation~\eqref{primaAtau} with respect to time
over the interval $(t,T)$, where $t$ is arbitrary in~$[0,T)$,
and use the definition~\eqref{backconv}.
Then, we account for the final condition~\eqref{cauchyAtau} and rearrange.
We then obtain the identity 
\Bsist
  && \iO (\pt+\tau\qt) v
  + \iG (\pGt+\tau\qGt) \vG
  + \iO \nabla(1*\qt) \cdot \nabla v
  + \iG \nablaG(1*\qGt) \cdot \nablaG\vG
  \non
  \\
  && = - \iO (1*(\psit\qt)) v
  - \iG (1*(\psiGt\,\qGt)) \vG
  + \iO (1*(\ubt\cdot\nabla\pt)) \, v
  \non
  \\
  && \quad {}
  + \iO (1*\juerg{\phi_3^\tau}) \, v
  + \iG (1*\juerg{\phi_4^\tau}) \, \vG
  + \iO \juerg{\phi_5^\tau} v
  + \iG \juerg{\phi_6^\tau} \vG \,.
  \label{intprimaAtau}
\Esist
This equality holds at any time and for every $(v,\vG)\in\calV$.
In the sequel, we use the notations
\Beq
  Q^t := \Omega \times (t,T)
  \aand
  \Sigma^t := \Gamma \times (t,T)
  \label{defbackQS}
\Eeq
for $t\in[0,T)$.

\def\intQt{\int_{Q^t}}
\def\intSt{\int_{\Sigma^t}}

\step
Basic estimate

For a future use, we notice the following property of the convolution \eqref{backconv}:
if $a\in L^\infty(Q)$, $b\in L^2(Q)$, $a_\Gamma\in L^\infty(\Sigma)$ and $b_\Gamma\in L^2(\Sigma)$, 
\juerg{then we have, for every $t\in[0,T)$,}
\Bsist
  && \intQt |1*(ab)|^2
  \leq T \, \norma a_\infty^2 \int_t^T \Bigl( \int_{Q^s} |b|^2 \Bigr) \, ds\,,
  \label{propconv}
  \\
  \separa
  && \intSt |1*(a_\Gamma b_\Gamma)|^2
  \leq T \, \norma{a_\Gamma}_\infty^2 \int_t^T \Bigl( \int_{S^s} |b_\Gamma|^2 \Bigr) \, ds\, .
  \label{propconvG}
\Esist
We write \eqref{intprimaAtau}, \eqref{secondaAtau}, and \eqref{primaAtau} at the time~$s$ 
and test them~by 
\Beq
  (\qt,\qGt)(s), \quad
  (\pt,\pGt)(s)-(\qt,\qGt)(s)
  \aand
  (\pt+\tau\qt,\pGt+\tau\qGt)(s),
  \non
\Eeq
respectively.
Then, we integrate over $(t,T)$ with respect to~$s$,
where $t\in[0,T)$ is arbitrary, obtaining the identities
\Bsist
  && \intQt (\pt+\tau\qt) \qt
  + \intSt (\pGt+\tau\qGt)\, \qGt
  + \intQt \nabla(1*\qt) \cdot \nabla\qt
  + \intSt \nablaG(1*\qGt) \cdot \nablaG\qGt
  \non
  \\
  && = - \intQt (1*(\psit\qt)) \qt
  - \intSt (1*(\psiGt\,\qGt))\, \qGt
  + \intQt (1*(\ubt\cdot\nabla\pt)) \, \qt
  \non
  \\
  && \quad {}
  + \intQt (1*\juerg{\phi_3^\tau}) \, \qt
  + \intSt (1*\juerg{\phi_4^\tau}) \, \qGt
  + \intQt \juerg{\phi_5^\tau} \, \qt
  + \intSt \juerg{\phi_6^\tau} \, \qGt\,,
  \non
  \\[2mm]
  \separa
  && \intQt \nabla\pt \cdot \bigl( \nabla\pt - \nabla\qt \bigr)
  + \intSt \nablaG\pGt \cdot \bigl( \juerg{\nablaG}\pGt - \juerg{\nablaG}\qGt \bigr)
  \non
  \\
  && = \intQt \qt (\pt - \qt)
  + \intSt \qGt \,(\pGt - \qGt) \,,
  \non
  \\[2mm]
  \separa
  && - \intQt \dt(\pt+\tau\qt) \, (\pt+\tau\qt) 
  - \intSt \dt(\pGt+\tau\qGt) \, (\pGt+\tau\qGt)
  \non
  \\
  && \quad {}
  + \intQt \bigl( \nabla\qt \cdot (\nabla\pt+\tau\nabla\qt \bigr) \bigr)
  + \intSt \bigl( \nablaG\qGt \cdot (\juerg{\nablaG}\pGt+\tau\nablaG\qGt \bigr) \bigr)
  \non
  \\
  && \quad {}
  + \intQt \psit\qt (\pt+\tau\qt)
  + \intSt \psiGt\,\qGt\, (\pGt+\tau\qGt)
  - \intQt \ubt \cdot \nabla\pt \, (\pt+\tau\qt)
  \non
  \\
  && = \intQt \juerg{\phi_3^\tau} \, (\pt+\tau\qt)
  + \intSt \juerg{\phi_4^\tau} \, (\pGt+\tau\qGt) .
  \non
\Esist
Now, we add these equalities to each other and account for \eqref{cauchyAtau}.
Since several terms cancel out, we obtain
\Bsist
  && \tau \intQt |\qt|^2
  + \tau \intSt |\qGt|^2
  + \frac 12 \iO |\nabla(1*\qt)(t)|^2
  + \frac 12 \iG |\nablaG(1*\qGt)(t)|^2
  \non
  \\
  && \quad {}
  + \intQt |\nabla\pt|^2
  + \intSt |\nablaG\pGt|^2
  + \frac 12 \iO |(\pt+\tau\qt)(t)|^2
  + \frac 12 \iG |(\pGt+\tau\qGt)(t)|^2
  \non
  \\
  && \quad {}
  + \intQt |\qt|^2
  + \intSt |\qGt|^2
  + \tau \intQt |\nabla\qt|^2 
  + \tau \intSt |\nablaG\qGt|^2 
  \non
  \\[1mm]
  \separa
  && = \frac 12 \iO |\juerg{\phi_5^\tau}|^2
  + \frac 12 \iG |\juerg{\phi_6^\tau}|^2
  \non
  \\
  && \quad {}
  - \intQt (1*(\psit\qt))\, \qt
  - \intSt (1*(\psiGt\,\qGt))\, \qGt
  + \intQt (1*(\ubt\cdot\nabla\pt)) \, \qt
  \non
  \\
  && \quad {}
  + \intQt (1*\juerg{\phi_3^\tau}) \, \qt
  + \intSt (1*\juerg{\phi_4^\tau}) \, \qGt
  + \intQt \juerg{\phi_5^\tau} \, \qt
  + \intSt \juerg{\phi_6^\tau} \, \qGt
  \non
  \\
  && \quad {}
  - \intQt \psit\qt (\pt+\tau\qt)
  - \intSt \psiGt\,\qGt \,(\pGt+\tau\qGt)
  + \intQt \ubt \cdot \nabla\pt \, (\pt+\tau\qt)
  \non
  \\
  && \quad {}
  + \intQt \juerg{\phi_3^\tau} \, (\pt+\tau\qt)
  + \intSt \juerg{\phi_4^\tau} \, (\pGt+\tau\qGt) .
  \non
\Esist
Every term on the \lhs\ is nonnegative. 
By simply using Young's inequality 
and \accorpa{propconv}{propconvG},
we see that the \rhs\ is estimated from above~by
\Bsist
  && \frac 12 \intQt |\qt|^2
  + \frac 12 \intSt |\qGt|^2
  + \frac 12 \intQt |\nabla\pt|^2
  \non
  \\
  && \quad {}
  + c \, \norma\psit_\infty^2 \intQt |1*\qt|^2
  + c \, \norma\psiGt_\infty^2 \intSt |1*\qGt|^2
  + c \, \norma\ubt_\infty^2 \intQt |1*\nabla\pt|^2
  \non
  \\
  && \quad {}
  + c \, \bigl( \norma\psit_\infty^2 + \norma\ubt_\infty^2 + \|\juerg{\phi_3^\tau}\|_2^2 \bigr) \intQt |\pt+\tau\qt|^2 
  \non
  \\
  && \quad {}
  + c \, \bigl( \norma\psiGt_\infty^2 + \|\juerg{\phi_4^\tau}\|_2^2 \bigr) \intSt |\pGt+\tau\qGt|^2 
  + c \, \sum_{i=3}^6 \|\juerg{\phi_i^\tau}\|_2^2\,,
  \non
\Esist
since $\norma{1*a}_2\leq c\,\norma a_2$ for every $L^2$ function~$a$.
By recalling \juerg{the uniform boundedness properties with respect to $\tau\in (0,1)$ of the controls $\ubt$ and of the functions defined in \eqref{defpsit} and \eqref{defphit}, 
we conclude from Gronwall's lemma that}
\Bsist
  && \norma{(\nabla\pt,\nablaG\pGt)}_{(\L2\calH)^3}
  + \norma{(\qt,\qGt)}_{\L2\calH}
  \non
  \\[1mm]
  && 
  + \,\norma{1*(\nabla\qt,\nablaG\qGt)}_{(\L\infty\calH)^3}
  + \norma{\pt+\tau\qt}_{\L\infty\calH}
  \leq c \,.
  \label{basic}
\Esist
Moreover, as \,\,$\norma{(\pt,\pGt)}_{\L2H}\leq c\,(\norma{\pt+\tau\qt}_{\L\infty\calH}+\norma{(\qt,\qGt)}_{\L2\calH})$,
we also have that
\Beq
  \norma{(\pt,\pGt)}_{\L2\calV} \leq c \,.
  \label{basicbis}
\Eeq

\vspace{5mm}\noindent
{\em Proof of Theorem~\ref{Necessary}:} \quad\,
\juerg{Let an optimal control $\ub\in\Uad$ for the original control problem with $\tau=0$ be given and $(\bar\mu,\bar\mu_\Gamma,
\bar\rho,\bar\rho_\Gamma)={\cal S}_0(\bar u)$ be the associated state.
By virtue of Theorem 4.2, we can pick a sequence $\{\tau_n\}\subset (0,1)$
with $\tau_n\searrow 0$ and, for any $n\in\enne$, an optimal control $u^{\tau_n}\in\Uad$, with associated
solution $(\mutn,\muGtn,\rhotn,\rhoGtn)={\cal S}_{\tau_n}(\utn)$ to the $\tau_n$-system \accorpa{primatau}{cauchytau}
for the adapted control problem $(\widetilde{\cal P}_{\tau_n})$ related to the cost functional~\eqref{cost2},
such that the following convergence properties hold true (see \eqref{conmu},\eqref{conrho}, and \eqref{strongconv}):}
\begin{align*}
&\juerg{\utn \to \ub
  \quad \hbox{strongly in $(\LQ2)^3$},}\\[1mm]
&\juerg{(\mutn,\muGtn) \to (\bar\mu,\bar\mu_\Gamma) \quad \mbox{weakly-star in }\,L^2(0,T;\calV),}\\[1mm]
&\juerg{(\rhotn,\rhoGtn) \to (\bar\rho,\bar\rho_\Gamma) \quad \mbox{weakly-star in } 
W^{1,\infty}(0,T;{\cal V}^*)\cap L^\infty(0,T;\calW)}\\
&\juerg{\hspace*{38.5mm} \mbox{and strongly in \,$C^0(\overline Q)\times C^0(\overline\Sigma)$}.}  
\end{align*} 
\juerg{But then also}
\begin{align*}
&\juerg{\psi^{\tau_n}=f''(\rhotn)\,\to\,f''(\bar\rho)=\psi \quad\mbox{strongly in }\,C^0(\overline Q),}\\[1mm]
&\juerg{\psi_\Gamma^{\tau_n}=\fG''(\rhoGtn)\,\to\,\fG''(\bar\rho_\Gamma)=\psi_\Gamma \quad
\mbox{strongly in }\,C^0(\overline\Sigma)\,,}
\end{align*} 
\juerg{and, likewise (recall \eqref{defphit}),}
\begin{align*}
&\juerg{\phi_3^{\tau_n}\to\phi_3\quad\mbox{strongly in } C^0(\overline Q),\qquad \phi_4^{\tau_n}\to\phi_4
\quad\mbox{strongly in } C^0(\overline\Sigma),}\\[1mm]
&\juerg{\phi_5^{\tau_n}\to \phi_5 \quad\mbox{strongly in } C^0(\overline\Omega), \qquad
\phi_6^{\tau_n}\to\phi_6 \quad\mbox{strongly in } C^0(\Gamma).}
\end{align*}
\juerg{Now observe that the dependence of the adapted and the non-adapted cost functionals on the state variables is 
the same. Therefore, nothing changes in the construction of the corresponding adjoint system,
which are both  given by \accorpa{primaAtau}{cauchyAtau}. In particular, the basic estimates
\eqref{basic} and \eqref{basicbis} are valid for the adjoint variables, and we thus may assume 
without loss of generality that there are $(p,p_\Gamma,q,q_\Gamma)$ such that}
\Bsist
  && (\ptn,\pGtn) \to (p,\pG)
  \quad \hbox{weakly in $\L2\calV$},
  \non
  \\[1mm]
  && (\qtn,\qGtn) \to (q,\qG)
  \quad \hbox{weakly in $\L2\calH$},
  \non
  \\[1mm]
  && 1*(\nabla\qtn,\nablaG\qGtn) \to 1*(\nabla q,\nablaG\qG)
  \quad \hbox{weakly-star in $(\L\infty\calH)^3$}.
  \non
\Esist
\juerg{We thus have the convergence properties}
\begin{align*}
&\juerg{\psi^{\tau_n}\,\qtn\,\to\,\psi\, q\quad\mbox{weakly in } L^2(Q), \qquad
\psi_\Gamma^{\tau_n}\,\qGtn\,\to\,\psiG\,\qG\quad\mbox{weakly in } L^2(\Sigma),}\\[1mm]
&\juerg {1*(\utn\cdot\nabla\ptn)\,\to\, 1*(\bar u\cdot\nabla p) \quad\mbox{weakly in } L^1(Q).}
\end{align*}
\juerg{Therefore, we may pass to the limit as $n\to\infty$ in the equations 
\eqref{intprimaAtau} and \eqref{secondaAtau}, written for $\tau=\tau_n$, $n\in\enne$,
to conclude that the quadruplet $(p,\pG,q,\qG)$ has the regularity properties
\eqref{regp}--\eqref{regq} and solves the system \eqref{primaA}--\eqref{secondaA}.} 

\juerg{It remains to show the validity of the variational inequality \eqref{necessary}.
In this regard, we observe that the variational inequalities for $({\cal P}_{\tau_n}\gianni)$
and $(\widetilde {\cal P}_{\tau_n})$ differ. Namely, for the adapted control problem it 
takes the form}
\Beq
  \intQ \bigl( \rho^{\tau_n} \, \nabla\ptn + \beta_7 \,\utn + (\utn - \ub) \bigr) \cdot 
	(v-\utn)
  \geq 0
  \quad \hbox{for every $v\in\Uad\,$}.
  \label{necessarytrick}
\Eeq
\juerg{Therefore, passage to the limit as $n\to\infty$ in \eqref{necessarytrick},
invoking the above convergence properties, yields the validity of \eqref{necessary}.
This concludes the proof of the assertion. \QED}

\Brem
\label{RemuniqA}
Coming back to \eqref{basic},
it is clear that the constant $c$ is proportional 
to the sum $\sum_{i=3}^6\juerg{\norma{\phi_i^\tau}_2}$
through a constant that depends only on the structure of the state system,
the control box and the initial datum.
In particular, if each $\juerg{\phi_i^\tau}$ vanishes, the inequality implies that
its \lhs\ vanishes as well.
Since the problem is linear, this is a uniqueness result for the inhomogeneous problem
associated to generic~$\phi_i$'s.
Therefore, if the same procedure were correct for $\tau=0$,
we would obtain a uniqueness result for the corresponding adjoint problem.
Coming back to the beginning of the procedure that led to~\eqref{basic},
we see that what is needed is the validity of the original equation 
\eqref{primaAtau} with $\tau=0$,
in addition to \accorpa{primaA}{secondaA}.
On the other hand, this is a consequence of \eqref{primaA} 
provided that $(q,\qG)\in\L2\calV$.
Indeed, in this case one is allowed to differentiate~\eqref{primaA}.
In conclusion, we \juerg{then would}  have a uniqueness result 
for the solution $(p,\pG,q,\qG)$ to the adjoint problem \accorpa{primaA}{secondaA}
whose component $(q,\qG)$ is smoother than required.
\Erem


\section*{Acknowledgments}
GG gratefully acknowledges some financial support from 
the GNAMPA (Gruppo Nazio\-nale per l'Analisi Matematica, 
la Probabilit\`a e le loro Applicazioni) of INdAM (Isti\-tuto 
Nazio\-nale di Alta Matematica) and the IMATI -- C.N.R.\ of Pavia.


\vspace{3truemm}


\vspace{3truemm}

\Begin{thebibliography}{10}

\bibitem{AbWi}
H. Abels, M. Wilke,
Convergence to equilibrium for the Cahn--Hilliard equation with a logarithmic free energy,
{\em Nonlinear Anal.} {\bf 67} (2007), 3176-3193.

\bibitem{Barbu}
V. Barbu,
``Nonlinear Differential Equations of Monotone Type in Banach Spaces'',
Springer,
London, New York, 2010.

\bibitem{Brezis}
H. Brezis,
``Op\'erateurs maximaux monotones et semi-groupes de contractions
dans les espaces de Hilbert'',
North-Holland Math. Stud.
{\bf 5},North-Holland,
Amsterdam,
1973. 

\bibitem{CaCo} 
L. Calatroni, P. Colli,
Global solution to the Allen--Cahn equation with singular potentials and
dynamic boundary conditions, {\it Nonlinear Anal.\/} {\bf 79} (2013), 12-27.

\bibitem{ChFaPr}
R. Chill, E. Fa\v sangov\'a, J. Pr\"uss, 
Convergence to steady state of solutions of the Cahn--Hilliard and Caginalp equations with dynamic boundary conditions, 
{\em Math. Nachr.} {\bf 279} (2006), 448-1462.

\bibitem{CFGS1}
P. Colli, M.~H. Farshbaf-Shaker, G. Gilardi, J. Sprekels, Optimal boundary control
of a viscous Cahn--Hilliard system with dynamic boundary condition and double
obstacle potentials, {\em SIAM J. Control Optim.} {\bf 53} (2015),  2696-2721.  

\bibitem{CFGS2}
P. Colli, M.~H. Farshbaf-Shaker, G. Gilardi, J. Sprekels, Second-order analysis of a 
boundary control problem for the viscous Cahn--Hilliard equation with dynamic
boundary conditions, {\em Ann. Acad. Rom. Sci. Math. Appl.} {\bf 7} (2015),
41-66.

\bibitem{CGLN}
P. Colli, G. Gilardi, P. Lauren\c cot, A. Novick-Cohen,
Uniqueness and long-time behaviour for the conserved phase-field system with memory,
{\em Discrete Contin. Dynam. Systems} {\bf 5} (1999), 375-390.

\bibitem{CGPS1}
P. Colli, G. Gilardi, P. Podio-Guidugli, J. Sprekels, 
Distributed optimal control of a nonstandard system of phase field equations,
{\em Contin. Mech. Thermodyn.} {\bf 24} (2012), 437-459.

\bibitem{CGPS2}
P. Colli, G. Gilardi, P. Podio-Guidugli, J. Sprekels,
Well-posedness and long-time behaviour for a nonstandard viscous Cahn--Hilliard system,
{\em SIAM J. Appl. Math.} {\bf 71} (2011), 1849-1870.

\bibitem{CGRS}
P. Colli, G. Gilardi, E. Rocca, J. Sprekels, Optimal distributed control of a
diffuse interface model of tumor growth, {\em Nonlinearity} {\bf 30} (2017), 2518-2546.

\bibitem{CGS1}
P. Colli, G. Gilardi, J. Sprekels, Analysis and optimal boundary control of a nonstandard 
system of phase field equations, {\em Milan J. Math.} {\bf  80} (2012), 119-149.

\bibitem{CGS8}
P. Colli, G. Gilardi, J. Sprekels, A boundary control problem for the viscous Cahn--Hilliard 
equation with dynamic boundary conditions, {\em Appl. Math. Optim.} {\bf 72} (2016), 195-225.

\bibitem{CGS9}
P. Colli, G. Gilardi, J. Sprekels, A boundary control problem for the pure 
Cahn--Hilliard equation with dynamic boundary conditions, {\em Adv. Nonlinear Anal.} 
{\bf 4} (2015), 311-325.

\bibitem{CGS10}
P. Colli, G. Gilardi, J. Sprekels, 
Global existence for a nonstandard viscous Cahn--Hilliard system with dynamic boundary condition, {\em SIAM J. Math. Anal.}
{\bf 49} (2017) 1732-1760.

\bibitem{CGS11}
P. Colli, G. Gilardi, J. Sprekels, Distributed optimal control of a nonstandard nonlocal phase field system, 
{\em AIMS Math.} {\bf 1} (2016), 246-281.

\bibitem{CGS12}
P. Colli, G. Gilardi, J. Sprekels, Distributed optimal control of a nonstandard nonlocal 
phase field system with double obstacle potential, {\em Evol. Equ. Control Theory} {\bf 6}
(2017), 35-58.

\bibitem{CGS13}
P. Colli, G. Gilardi, J. Sprekels,
On a Cahn--Hilliard system with convection and
dynamic boundary conditions, \giorgio{{\em Ann. Mat. Pura. Appl. (4)}, Online First 2018,
DOI 10.1007/s10231-018-0732-1.}

\bibitem{CGS14}
P. Colli, G. Gilardi, J. Sprekels,
Optimal velocity control of a viscous Cahn--Hilliard system 
with convection and dynamic boundary conditions, \giorgio{{\em SIAM J. Control. Optim.}, to appear 2018}.

\bibitem{CGS16}
P. Colli, G. Gilardi, J. Sprekels,
Optimal velocity control of a convective Cahn--Hilliard system 
with double obstacles and dynamic boundary conditions: a `deep quench' approach, 
\giorgio{Preprint arxiv:1709.03892 [math.AP] (2017), pp. 1-30.}

\bibitem{CGS17}
P. Colli, G. Gilardi, J. Sprekels,
On the \longtime\ \bhv\ of a viscous Cahn--Hilliard system
with convection and dynamic boundary conditions, 
Preprint arxiv:1803.04318 [math.AP] (2018), pp. 1-20.

\bibitem{CS}
P. Colli, J. Sprekels, Optimal boundary control of a nonstandard Cahn--Hilliard system 
with dynamic boundary condition and double obstacle inclusions, \giorgio{in: ``Solvability, Regularity, and Optimal Control of Boundary 
Value Problems for PDEs'' (P. Colli, A. Favini, E. Rocca, G. Schimperna, J. Sprekels, eds.), pp. 151-182, Springer INdAM Series {\bf 22},
Springer, Milan 2017.}

\giorgio{
\bibitem{Duan}
N. Duan, X. Zhao, Optimal control for the multi-dimensional viscous Cahn--Hilliard equation, {\em Electronic J. Differ. Equations 2015},
Paper No. 165, 13 pp.
}

\bibitem{EfMiZe}
M. Efendiev, A. Miranville, S. Zelik, 
Exponential attractors for a singularly perturbed Cahn--Hilliard system. 
{\em Math. Nachr.} {\bf 272} (2004), 11-31.

\bibitem{EfGaZe}
M. Efendiev, H. Gajewski, S. Zelik, 
The finite dimensional attractor for a 4th order system of the Cahn--Hilliard type with a supercritical nonlinearity,
{\em Adv. Differential Equations} {\bf 7} (2002), 1073-1100.

\bibitem{FRS}
S. Frigeri, E. Rocca, J. Sprekels, Optimal distributed control of a nonlocal 
Cahn--Hilliard/Navier--Stokes system in two dimensions.
{\em SIAM J. Control. Optim.} {\bf 54} (2016), 221-250.

\giorgio{
\bibitem{Fukao}
T. Fukao, N. Yamazaki, A boundary control problem for the equation and dynamic boundary condition of Cahn--Hilliard type,
in: ``Solvability, Regularity, and Optimal Control of Boundary Value Problems for PDEs'' (P. Colli, A. Favini, E. Rocca, 
G. Schimperna, J. Sprekels, eds.), pp. 255-280, Springer INdAM Series {\bf 22},
Springer, Milan 2017.} 

\bibitem{Gal1}
C.\,G. Gal,
Exponential attractors for a Cahn--Hilliard model in bounded domains with permeable walls,
{\em Electron. J. Differential Equations} {\bf 2006} (2006), 1-23.

\bibitem{Gal2}
C. Gal,
Well-posedness and long time behavior of the non-isothermal viscous Cahn--Hilliard equation with dynamic boundary conditions,
{\em Dyn. Partial Differ. Equ.} {\bf 5} (2008), 39-67.

\bibitem{GaGr1}
C.\,G. Gal, M. Grasselli,
Asymptotic behavior of a Cahn--Hilliard--Navier--Stokes system in~2D,
{\em Ann. Inst. H. Poincar\'e Anal. Non Lin\'eaire } {\bf 27} (2010), \giorgio{401-436}.

\bibitem{GaGr2}
C.\,G. Gal, M. Grasselli, 
Instability of two-phase flows: A lower bound on the dimension of the global attractor 
of the Cahn--Hilliard--Navier--Stokes system,
{\em Physica~D: Nonlinear Phenomena}
{\bf 240} (2011), 629-635.

\bibitem{GiMiSchi} 
G. Gilardi, A. Miranville, and G. Schimperna,
On the Cahn--Hilliard equation with irregular potentials and dynamic boundary conditions,
{\it Commun. Pure Appl. Anal.\/} 
{\bf 8} (2009), 881-912.

\bibitem{GiRo}
G. Gilardi, E. Rocca,
Well posedness and long time behaviour for a singular phase field system of conserved type
IMA J. Appl. Math. {\bf 72} (2007) 498-530.

\bibitem{GrPeSch}
M. Grasselli, H. Petzeltov\'a, G. Schimperna,  
Asymptotic  behavior  of  a  nonisothermal  viscous  Cahn--Hilliard  equation  with  inertial  term,
{\em J. Differential Equations} {\bf 239} (2007), 38-60.

\bibitem {HHKK}
M. Hinterm\"uller, M. Hinze, C. Kahle, T. Keil, A goal-oriented dual-weighted adaptive 
finite element approach for the optimal control of a nonsmooth Cahn--Hilliard--Navier--Stokes system,
WIAS Preprint No. 2311, Berlin, 2016.

\bibitem{HiKeWe}
M. Hinterm\"uller, T. Keil, D. Wegner, Optimal control of a semidiscrete 
Cahn--Hilliard--Navier--Stokes  system  with  non-matched  fluid  densities, \giorgio{
{\em SIAM J. Control Optim.} {\bf 55} (2018), 1954-1989.}

\bibitem{HW1}
M. Hinterm\"uller, D. Wegner,
Distributed optimal control of the Cahn--Hilliard system including the case of a double-obstacle homogeneous free energy density, 
{\em SIAM J. Control Optim.} {\bf 50}
(2012), 388-418.

\bibitem{HW2}
M. Hinterm\"uller, D. Wegner,
Optimal control of a semidiscrete Cahn--Hilliard--Navier--Stokes system, 
{\em SIAM J. Control Optim.} {\bf 52} (2014), 747--772.

\bibitem{HW3}
 M. Hinterm\"uller, D. Wegner, Distributed and boundary control problems for the semidiscrete 
Cahn--Hilliard/Navier--Stokes system with nonsmooth Ginzburg--Landau energies. \giorgio{{\em Topological Optimization and Optimal Transport,
Radon Series on Computational and Applied Mathematics} {\bf 17} (2017), 40-63.}

\bibitem{JiWuZh}
J. Jiang, H. Wu, S. Zheng,
Well-posedness and long-time behavior of a non-autonomous Cahn--Hilliard--Darcy system with mass source modeling tumor growth
{\em J. Differential Equations} {\bf 259} (2015), 3032-3077.

\bibitem{LiZho}
D. Li, C. Zhong,
Global attractor for the Cahn--Hilliard system with fast growing nonlinearity,
{\em J. Differential Equations} {\bf 149} (1998), 191-210.

\giorgio{
\bibitem{Medjo}
T. Tachim Medjo, Optimal control of a Cahn--Hilliard--Navier--Stokes model with state constraints, {\em J. Convex Anal.} {\bf 22} (2015),
1135-1172.
}

\bibitem{Mir1}
A. Miranville, 
Asymptotic behavior of the Cahn--Hilliard--Oono equation, 
{\em J. Appl. Anal. Comput.} {\bf 1} (2011), 523-536.

\bibitem{Mir2}
A. Miranville, 
Asymptotic behavior of a generalized Cahn--Hilliard equation with a proliferation term, 
{\em Appl. Anal.} {\bf 92} (2013), 1308-1321.

\bibitem{Mir3}
A. Miranville, 
Long-time behavior of some models of Cahn--Hilliard equations in deformable continua, 
{\em Nonlinear Anal.} {\bf 2} (2001), 273-304.

\bibitem{MiZe2}
A. Miranville, S. Zelik,
Robust exponential attractors for Cahn--Hilliard type equations with singular potentials,
{\em Math. Meth. Appl. Sci.} {\bf 27} (2004), 545-582.

\bibitem{MiZe1}
A. Miranville, S. Zelik, 
Exponential attractors for the Cahn--Hilliard equation with dynamic boundary conditions, 
{\em Math. Models Appl. Sci.} {\bf 28} (2005), 709-735.

\bibitem{PrVeZa}
J. Pr\"uss, V. Vergara, R. Zacher,
Well-posedness and long-time behaviour for the non-isothermal Cahn--Hilliard equation with memory,
{\em Discrete Contin. Dyn. Syst. Series~A} {\bf 26} (2010), 625-647.

\bibitem{RS}
E. Rocca, J. Sprekels, Optimal distributed control of a nonlocal convective Cahn--Hilliard equation by the 
velocity in three dimensions, {\em SIAM J. Control Optim.} {\bf 53} (2015), 1654-1680.

\bibitem{Seg}
A. Segatti, 
On the hyperbolic relaxation of the Cahn--Hilliard equation in~3D: approximation and long time behaviour,
{\em Math. Models Methods Appl. Sci.} {\bf 17} (2007), 411-437.

\bibitem{Simon}
J. Simon,
Compact sets in the space $L^p(0,T; B)$,
{\it Ann. Mat. Pura Appl.~(4)\/} 
{\bf 146}, (1987), 65-96.

\bibitem{WN}
Q.-F.  Wang,  S.-i.  Nakagiri, 
Weak  solutions  of  Cahn--Hilliard  equations  having  
forcing terms and optimal control problems, Mathematical models in functional equations
(Japanese) (Kyoto, 1999), Surikaisekikenkyusho Kokyuroku No. 1128 (2000), 172-180.

\bibitem{WaWu}
X.-M. Wang, H. Wu, 
Long-time behavior for the Hele--Shaw--Cahn--Hilliard system,
{\em Asymptot. Anal.} {\bf 78} (2012), 217-245.

\bibitem{WuZh}
H. Wu, S. Zheng, 
Convergence to equilibrium for the Cahn--Hilliard equation with dynamic boundary conditions, 
{\em J. Differential Equations} {\bf 204} (2004), 511-531.

\bibitem{ZhaoLiu}
X. Zhao, C. Liu,
On the existence of global attractor for 3D viscous Cahn--Hilliard equation,
{\em Acta Appl. Math.} {\bf 138} (2015), 199-212.

\bibitem{ZL1}
X. Zhao, C. Liu, 
Optimal control of the convective Cahn--Hilliard equation, 
{\em Appl. Anal.} {\bf 92} (2013), 1028-1045.

\bibitem{ZL2}
X. Zhao, C. Liu,
Optimal control for the convective Cahn--Hilliard equation in 2D case,
{\em Appl. Math. Optim.} {\bf 70} (2014), 61-82.

\giorgio{
\bibitem{JZheng1}
J. Zheng, Time optimal controls of the Cahn--Hilliard equation with internal control,
{\em Optimal Control Appl. Methods} {\bf 60} (2015), 566-582.
}

\giorgio{
\bibitem{ZW}
J. Zheng, Y. Wang, Optimal control problem for Cahn--Hilliard equations with state
constraints, {\em J. Dyn. Control Syst.} {\bf 21} (2015), 257-272.
}

\bibitem{Zheng}
S. Zheng,
Asymptotic behavior of solution to the Cahn--Hilliard equation,
{\em Appl. Anal.} {\bf 23} (1986), 165-184.

\End{thebibliography}

\End{document}
